\documentclass[a4paper,10pt]{amsart}

\usepackage[english]{babel}

\usepackage{verbatim}

\usepackage{amssymb} 

\usepackage{tikz, tikz-3dplot} 
\usepackage{tikz-cd}

\usepackage{lmodern} 

\usepackage{fancyvrb} 

\usepackage{url} 

\usepackage{hyperref} 

\newtheorem{thm}{Theorem}[section]
\newtheorem{trm}[thm]{Theorem}

\newtheorem{prop}[thm]{Proposition}
\newtheorem{lem}[thm]{Lemma}
\newtheorem{lma}[thm]{Lemma}

\theoremstyle{definition}

\newtheorem{defi}[thm]{Definition}

\newtheorem{ex}[thm]{Example}

\theoremstyle{remark}
\newtheorem{rem}[thm]{Remark}

\makeatletter
\let\c@equation\c@thm
\makeatother
\numberwithin{equation}{section}

\newcommand{\C}{\mathbb{C}}

\newcommand{\Q}{\mathbb{Q}}
\newcommand{\Z}{\mathbb{Z}}
\newcommand{\N}{\mathbb{N}}

\newcommand{\f}{\textbf}
\newcommand{\mr}{\mathrm}
\newcommand{\mb}{\mathbf}
\newcommand{\bdim}{\mathbf{dim}\,}
\newcommand{\fa}{ \ \mathrm{for}\ \mathrm{all}\ }

\title[Desingularizations of Cyclic Quiver Grassmannians]{Desingularizations of Quiver Grassmannians for the Equioriented Cycle Quiver}

\author{Alexander P\"{u}tz and Markus Reineke}

\begin{document}

\begin{abstract}
We construct torus equivariant desingularizations of quiver Grassmannians for arbitrary nilpotent representations of an equioriented cycle quiver. We apply this to the computation of their torus equivariant cohomology.
\end{abstract}

\maketitle
\section{Introduction}

Quiver Grassmannians are projective varieties parametrizing subrepresentations of quiver representations. Originating in the geometric study of quiver representations \cite{Schofield1992} and in cluster algebra theory \cite{CC2006}, they have been applied extensively in recent years in a Lie-theoretic context, namely as a fruitful source for degenerations of (affine) flag varieties \cite{CFFFR2017,CFR2013,FFR2017,Pue2020}. This approach allows for an application of homological methods from the representation theory of quivers to the study of such degenerate structures.

The resulting varieties being typically singular, a construction of natural desingularizations is very desirable. For quiver Grassmannians of representations of Dynkin quivers this was accomplished in \cite{CFR2013} building on \cite{FF13}, and for Grassmannians of subrepresentations of loop quivers in \cite{FFR2017}.

In the present paper, we synthesize the approaches of \cite{CFR2013,FFR2017} and construct desingularizations of quiver Grassmannians for nilpotent representations of equioriented cycle quivers, thereby, in particular, desingularizing degenerate affine flag varieties \cite{Pue2020}.

As an important application, this allows us to describe the equivariant cohomology of degenerate affine flag varieties and more general quiver Grassmannians, in continuation of the programme started in \cite{LaPu2020,LaPu2021}.


In the first section, we recall some background material on quiver Grassmannians for nilpotent representations of the equioriented cycle quiver. In Section~\ref{sec:desing} we give an explicit construction for desingularizations of quiver Grassmannians for nilpotent representations of the equioriented cycle, along the lines of \cite{CFR2013}. We prove that the desingularization has a particularly favourable geometric structure, namely it is isomorphic to a tower of Grassmann bundle. Consequently, it, admits a cellular decomposition which is compatible with the cellular decomposition of the singular quiver Grassmannian. In Section~\ref{sec:equivariant-cohomology}, we recall the definition of certain torus actions on cyclic quiver Grassmannians, together with the necessary framework to compute torus equivariant cohomology. 
Finally, in Section~\ref{sec:application} we prove that the desingularization is equivariant with respect to the torus action as introduced in \cite{LaPu2020}. This allows to use the construction from that paper for the computation of torus equivariant cohomology to all quiver Grassmannians for nilpotent representations of the equioriented cycle.

\section{Quiver Grassmannians for the Equioriented Cycle}\label{sec:quiver-grass-cycle}
In this section we recall some definitions concerning quiver Grassmannians and representations of the equioriented cycle. We refer to \cite{Kirillov2016,Schiffler2014} for general representation theoretic properties, and to \cite{CFR2012} for basic properties of quiver Grassmannians.

\subsection{Generalities on quiver representations}
Let $Q$ be a quiver, consisting of a set of vertices $Q_0$ and a set of arrows $Q_1$ between the vertices. A $Q$-representation $M$ consists of a tuple of $\C$-vector spaces $M^{(i)}$ for $i \in Q_0$ and tuple of linear maps $M_\alpha : M^{(i)} \to M^{(j)}$ for $(\alpha:i\to j)\in Q_1$. We denote the category of finite dimensional $Q$-representations by $\mr{rep}_\C(Q)$. The morphisms between two objects $M$ and $N$ are tuples of linear maps $\varphi_i : M^{(i)} \to N^{(i)}$ $i \in Q_0$ such that 
\( \varphi_j \circ M_\alpha = N_\alpha \circ \varphi_i\) holds for all $(\alpha:i\to j) \in Q_1$. 
\begin{defi}
For $M \in \mr{rep}_\C(Q)$ and $\mb{e} \in \N^{Q_0}$, the \f{quiver Grassmannian} $\mr{Gr}_\mb{e}(M)$ is the closed subvariety of $\prod_{i\in Q_0}{\rm Gr}_{e_i}(M^{(i)})$ of all subrepresentations $U$ of $M$ such that $\dim_\C U^{(i)}=e_i$ for $i \in Q_0$.
\end{defi}
For an isomorphim class $[N]$ of $Q$-representations, the stratum $\mathcal{S}_{[N]}$ is defined as the set of all points (that is, subrepresentations) $U \in \mr{Gr}_\mb{e}(M)$ such that $U$ is isomorphic to $N$. By \cite[Lemma 2.4]{CFR2012}, $\mathcal{S}_{[N]}$ is locally closed and irreducible. If there are only finitely many isomorphism classes of subrepresentations of $M$, as will be the case in the following, the $\mathcal{S}_{[N]}$ thus define a finite stratification of the quiver Grassmannians.

Every basis $B$ of $M \in \mr{rep}_\C(Q)$ consists of bases 
\[ B^{(i)} =  \big\{ v^{(i)}_k \ \big\vert \ k \in [m_i] \big\} \]
for each vector space $M^{(i)}$ of the $Q$-representation $M$, where $m_i := \dim_\C M^{(i)}$ for all $i \in Q_0$, and $[m]:=\{1,\dots,m\}$.
\begin{defi}
Let $M \in \mr{rep}_\C(Q)$ and $B$ a basis of $M$. The \f{coefficient quiver} $Q(M,B)$ consists of: 
\begin{itemize}
\item[(QM0)] the vertex set  $Q(M,B)_0=B$,
\item[(QM1)] the set of arrows $Q(M,B)_1$, containing $(\tilde{\alpha}: v_k^{(i)} \to v_\ell^{(j)})$ if and only if $(\alpha:i\to j) \in Q_1$ and the coefficient of $v_\ell^{(j)}$ in $M_\alpha v_k^{(i)}$ is non-zero.
\end{itemize}
\end{defi}

\subsection{Representations of the Equioriented Cycle}\label{sec:reps-cycle}
By $\Delta_n$ we denote the equioriented cycle quiver on $n$ vertices. Hence the set of arrows and the sets of vertices are in bijection with $\Z_n := \Z/n\Z$; more precisely, we have $(\Delta_n)_0=\Z/n\Z$ and arrows $i\rightarrow i+1$ for all $i\in \Z/n\Z$. Here and in the following, we consider all indices modulo $n$ unless specified differently.

A $\Delta_n$-representation $M$ is called \f{$N$-nilpotent} for $N \in \N$ if
\[ M_{\alpha_{i+N-1}} \circ  M_{\alpha_{i+N-2}} \circ \dots  \circ M_{\alpha_{i +1}} \circ M_{\alpha_i} = 0 \]
for all $i \in \Z_n$, i.e. all concatenations of the maps of $M$ along the arrows of $\Delta_n$ vanish after at most $N$ steps. $M$ is called \f{nilpotent} if it is $N$-nilpotent for some $N$.
\begin{ex}\label{ex:NilpoIndcpDeltan}Let $i\in \Z_n$ and let $\ell\in\Z_{\geq 1}$. Consider the $\C$-vector space $V$ with basis $B = \{b_1, \ldots, b_{\ell}\}$ equipped with the $\Z_n$-grading given by $\deg(b_k)=i+k-1 \in\Z_n$. Take the operator $A\in \mathrm{End}(V)$ uniquely determined by setting $Aw_{k}=w_{k+1}$ for any $k < \ell$ and $Aw_{\ell}=0$. It is immediate to check that the corresponding $\Delta_n$-representation is nilpotent. We denote this representation by $U_i(\ell)$.
\end{ex}
\begin{prop} (\cite[Theorem~7.6.(1)]{Kirillov2016})
Every indecomposable nilpotent representation of $\Delta_n$ is isomorphic to some $U_i(\ell)$. 
\end{prop}
\begin{ex}
Observe that the basis $B$ from Example~\ref{ex:NilpoIndcpDeltan} can be obviously rearranged into the union of ordered bases 
\[ B^{(i)} = \big\{ v^{(i)}_r : r \in [k_i] \big\} \quad \mr{for} \ i\in\Z_n, \]
where $k_i$ is the number of elements $b \in B$ with degree $\mr{deg}(b)=i$.
With respect to $B$, the coefficient quiver of $U_i(\ell)$ has the form:
\begin{center}
\begin{tikzpicture}[scale=.56]

  \node at ($(0,0)+(-306:4.7)$) {\tiny $v^{(i)}_1$};
  \node at ($(0,0)+(-133:1.6)$) {\tiny $v^{(j)}_{k_j}$};

\foreach \ang\dist in {-300/4.5, -330/4.4166, 0/4.3333, -30/4.25, -60/4.1666, -90/4.0833, -120/4.0, -150/3.9166, -180/3.8333, -210/3.75, -240/3.6666, -270/3.5833, -300/3.5, -330/3.4166, 0/3.3333, -30/3.25}{
  \draw[fill=black] ($(0,0)+(\ang:2+0.5*\dist)$) circle (.08);
}

\foreach \ang\dist in {-60/2.1666, -90/2.0833, -120/2.0, -150/1.9166, -210/1.75, -240/1.6666, -270/1.5833, -300/1.5, -330/1.4166, -0/1.3333, -30/1.25, -60/1.1666, -90/1.0833, -120/1.0, -150/0.9166, -210/0.75, -240/0.6666, -270/0.5833, -300/0.5, -330/0.4166, 0/0.3333, -30/0.25, -60/0.1666, -90/0.0833, -120/0.0}{
  \draw[fill=black] ($(0,0)+(\ang:2+0.5*\dist)$) circle (.08);
}

\foreach \ang\dist in {-300/4.5, 0/4.3333, -30/4.25, -60/4.1666, -90/4.0833, -120/4.0, -150/3.9166, -210/3.75, -240/3.6666, -270/3.5833, -300/3.5, 0/3.3333}{
  \draw[->,shorten <=1pt, shorten >=0pt] ($(0,0)+(\ang:2+0.5*\dist)$) arc (\ang-1.15:\ang-30:2+0.5*\dist);
}

\foreach \ang\dist in {-60/2.1666, -90/2.0833, -120/2.0, -150/1.9166, -210/1.75, -240/1.6666, -270/1.5833, -300/1.5, -0/1.3333, -30/1.25, -60/1.1666, -90/1.0833, -120/1.0, -150/0.9166 , -210/0.75, -240/0.6666, -270/0.5833, -300/0.5, 0/0.3333, -30/0.25, -60/0.1666, -90/0.0833}{
  \draw[->,shorten <=1pt, shorten >=0pt] ($(0,0)+(\ang:2+0.5*\dist)$) arc (\ang-1.15:\ang-29.5:2+0.5*\dist);
}

\foreach \ind in {1,2,3,4,5,6,7,8,9,10,11,12,13,14,15,16,17,18,19,20,21,22,23,24,25,26,27,28,29,30,31,32,33,34,35,36}{  
  \draw[fill=black] ($(0,0)+(-30-360*\ind/72:3.5433-\ind/144+11/144)$) circle (.015);
  \draw[fill=black] ($(0,0)+(-30-360*\ind/72-180:3.05-\ind/144+44/144)$) circle (.015);
}

\foreach \ind in {1,2,3,4,5,6}{  
  \draw[fill=black] ($(0,0)+(-30-360*\ind/72:3.05-\ind/144+11/144)$) circle (.015);
}

\foreach \ind in {2,3,4}{  
  \draw[fill=black] ($(0,0)+(180-360*\ind/72:3.629-\ind/144+44/144)$) circle (.015);
  \draw[fill=black] ($(0,0)+(180-360*\ind/72:3.629-1-\ind/144+44/144)$) circle (.015);
  \draw[fill=black] ($(0,0)+(180-360*\ind/72:3.629-1.5-\ind/144+44/144)$) circle (.015);
}

\foreach \ind in {2,3,4}{  
  \draw[fill=black] ($(0,0)+(30-360*\ind/72:4.129-\ind/144+11/144)$) circle (.015);
  \draw[fill=black] ($(0,0)+(30-360*\ind/72:4.129-0.5-\ind/144+11/144)$) circle (.015);
  \draw[fill=black] ($(0,0)+(30-360*\ind/72:4.129-1.5-\ind/144+11/144)$) circle (.015);
  \draw[fill=black] ($(0,0)+(30-360*\ind/72:4.129-2-\ind/144+11/144)$) circle (.015);

}

\end{tikzpicture}
\end{center}
\end{ex}
By  \cite[Theorem~1.11]{Kirillov2016}, every nilpotent $\Delta_n$-representation is isomorphic to a $\Delta_n$-representation of the form:
\[M := \bigoplus_{i\in \Z_n} \bigoplus_{\ell \in [N]} U_i(\ell) \otimes \C^{d_{i,\ell}},\]
with $d_{i,\ell} \in \Z_{\geq 0}$ for all $i \in \Z_n$ and $\ell \in [N] $. 

Let $\C\Delta_n$ be the path algebra of $\Delta_n$ and define the path 
\[ p_i(N):=(i+N \vert \alpha_{i+N-1} \circ \alpha_{i+N-2} \circ \dots \circ \alpha_{i+1} \circ\alpha_{i} \vert i )\]
for all $i \in \Z_n$ and some fixed $N \in \N$. We define the path algebra ideal  
\[\mr{I}_N := \langle \, p_i(N) : i \in \Z_n \, \rangle \subset \C \Delta_n, \]
generated by all paths of length $N$, and we denote the truncated path algebra $\C\Delta_n/\mr{I}_N$ by $A_n^{(N)}$. The following is a special case of \cite[Theorem~5.4]{Schiffler2014}.
\begin{prop}
The category $\mr{rep}_\C(\Delta_n,\mr{I}_N)$ of bounded quiver representations is equivalent to the category $\mr{mod}_\C(A_n^{(N)})$ of modules over the truncated path algebra. 
\end{prop}
\begin{rem}
The $U_i(N)$ for $i \in \Z_n$ are the longest indecomposable nilpotent representations in $\mr{rep}_\C(\Delta_n,\mr{I}_N)$.
\end{rem}
Let $P_i \in \mr{rep}_\C(\Delta_n,\mr{I}_N)$ be the $\Delta_n$-representation starting at vertex $i \in \Z_n$ corresponding to the projective indecomposable $A_n^{(N)}$-module, and similarly let $I_j \in \mr{rep}_\C(\Delta_n,\mr{I}_N)$ be the $\Delta_n$-representation ending at vertex $j \in \Z_n$ corresponding to the injective indecomposable $A_n^{(N)}$-representation. We can identify bounded projective and bounded injective representations of the cycle, via indecomposable nilpotent representations (c.f. \cite[Proposition~4.2]{Pue2020}):
\begin{prop}
\label{prop:proj_inj_dual}
For $n,N \in \N$ and all $i,j \in \Z_n$ the projective and injective representations $P_i$ and $I_j$ of the bound quiver $(\Delta_n,\mr{I}_N)$ satisfy 
\[ P_i \cong U_{i}(N) \cong I_{i+N-1} \quad \mr{and} \quad  I_j \cong U_{j-N+1}(N) \cong P_{j-N+1}.\]
\end{prop}
\begin{rem}
In particular, every indecomposable nilpotent  $\Delta_n$-representation $U_i(N)$ is projective and injective in $\mr{rep}_\C(\Delta_n,\mr{I}_N)$ . If we want to emphasize the injective nature of an indecomposable $\Delta_n$-representation we sometimes use the notation $U(j;\ell) := U_{j-\ell+1}(\ell)$.
\end{rem}
\subsection{Parametrization of Irreducible Components}
In Section~\ref{sec:desing} we will construct desingularizations of all quiver Grassmannians for nilpotent representations of the equioriented cycle, which requires knowledge of their irreducible components. Let us first recall the approach: since there are only finitely many isomorphism classes of nilpotent $\Delta_n$-representations in any fixed dimension, the stratification of every quiver Grassmannian into strata $\mathcal{S}_{[N]}$ is finite. Since the strata are irreducible, the irreducible components of quiver Grassmannians are therefore of the form $\overline{\mathcal{S}_{[N]}}$ for certain isomorphism classes $[N]$, which provide a natural labelling (and a canonical representative) of the components.

For arbitrary nilpotent representations of the equioriented cycle the structure of the irreducible components 
is not known.
In the special case that all indecomposable direct summands of the $\Delta_n$-representation $M$ have length  $N= \omega n$ and $\mb{e} = (\omega k, \dots , \omega k ) \in \Z^n$, we have an explicit description of the irreducible components of the quiver Grassmannian $\mr{Gr}_\mb{e} (M)$ \cite[Lemma~4.10]{Pue2020}:
\begin{lma}
\label{lma:irr_comp}
Let $M$ denote the $\Delta_n$-representation $\oplus_{i \in \Z_n} U(i;\omega n) \otimes \C^{d_i}$ with $d_i \in \Z_{\geq 0}$ for all $i \in \Z_n$, define $m :=\sum_{i \in \Z_n} d_i$ and $\mb{e} := (\omega k, \dots , \omega k ) \in \Z^n$. The irreducible components of $\mr{Gr}_\mb{e}(M)$ are in bijection with the set
\[ C_k(\mb{d}) := \Big\{ \mb{p} \in \Z_{\geq 0}^{n} : p_i \leq d_i \ \mr{for} \ \mr{all} \ i \in \Z_n, \sum_{i \in \Z_n} p_i = k \Big\} \]
and they all have dimension $\omega k(m-k)$. 
\end{lma}
\begin{rem}
A representative of the open dense stratum in the irreducible component corresponding to $\mb{p} \in C_k(\mb{d})$ is  
\[ U_\mb{p} := \bigoplus_{i \in \Z_n} U(i;\omega n)\otimes \C^{p_i}.\] 
\end{rem}
\begin{ex}\label{ex:generic-subrep-special-case}
Let $d_i=1$ for all $i \in \Z_n$. Then by Lemma~\ref{lma:irr_comp} the irreducible components are parametrized by the $k$-element subsets of $[n]$ and the representatives are 
\[ \bigoplus_{j \in I} U(j;\omega n)\]
for $I \in \binom{[n]}{k}$. The dimension of the irreducible components is $\omega k(n-k)$. 
\end{ex}
\section{Construction of the Desingularization}\label{sec:desing}
The approach to the construction of desingularizations of quiver Grassmannians for the equioriented cycle quiver carried out in this section is a synthesis of the approach of \cite{CFR2013} for Dynkin quivers and the approach of \cite{FFR2017} for the loop quiver. We will construct another quiver for which certain quiver Grassmannians yield desingularizations, which relies on certain favourable homological properties similar to those in \cite[Section~4]{CFR2013}.
\subsection{Bounded Representations of the Equioriented Cylinder}
In this subsection we introduce a map $ \Lambda : \mr{rep}_\C(\Delta_n,\mr{I}_N) \to \mr{rep}_\C(Q,\mr{I})$ for some bound quiver $(Q,\mr{I})$ such  that each quiver Grassmannian associated to $\Lambda(M)$ is smooth for all $M \in \mr{rep}_\C(\Delta_n,\mr{I}_N)$. We start with the definition of $Q$ and the ideal $\mr{I}$.
Let $\hat{\Delta}_{n,N}$ be the quiver with vertices $\big(\hat{\Delta}_{n,N}\big)_0 =\{(i,k) \, : \; i \in \Z_n, \ k \in [N] \}$ and arrows  
\begin{align*}
\big(\hat{\Delta}_{n,N}\big)_1 = &\big\{ \alpha_{i,k} : (i,k) \to (i,k+1) \ : \ i \in \Z_n, \ k \in [N-1] \big\}\cup\\
  &\big\{ \beta_{i,k}: (i,k) \to (i+1,k-1) \ : \ i \in \Z_n, \ k \in [N] \setminus \{ 1\} \big\}
 \end{align*}
which we call an equioriented cylinder quiver. We define $\hat{\mr{I}}_{n,N}$ as the ideal in the path algebra $\C \hat{\Delta}_{n,N}$ generated by the relations
 \[ \beta_{i,k+1} \circ \alpha_{i,k} \equiv \alpha_{i+1,k-1} \circ \beta_{i,k} \quad \mr{and} \quad \alpha_{i+1,N-1} \circ \beta_{i,N} \equiv 0 \]
 for all $i \in \Z_n$ and all $k \in [N-1]\setminus \{ 1\}$.

 \begin{ex}
$\hat{\Delta}_{4,4}$ is the following quiver: 
\begin{center}
\begin{tikzpicture}[scale=1.6]
\node at (0,1) {$(1,1)$};\node at (0.5,0.75) {$\alpha_{1,1}$}; 
\node at (1,0) {$(2,1)$};\node at (0.7,-0.45) {$\alpha_{2,1}$};
\node at (0,-1) {$(3,1)$};\node at (-0.5,-0.75) {$\alpha_{3,1}$};
\node at (-1,0) {$(4,1)$};\node at (-0.7,0.45) {$\alpha_{4,1}$};

\draw[arrows={-angle 90}, shorten >=15, shorten <=15]  (0,1) -- (1,1);
\draw[arrows={-angle 90}, shorten >=15, shorten <=15]  (1,0) -- (1,-1);
\draw[arrows={-angle 90}, shorten >=15, shorten <=15]  (0,-1) -- (-1,-1);
\draw[arrows={-angle 90}, shorten >=15, shorten <=15]  (-1,0) -- (-1,1);

\node at (1,1) {$(1,2)$};\node at (0.7,0.45) {$\beta_{1,2}$};
\node at (1,-1) {$(2,2)$};\node at (0.5,-0.75) {$\beta_{2,2}$};
\node at (-1,-1) {$(3,2)$};\node at (-0.7,-0.45) {$\beta_{3,2}$};
\node at (-1,1) {$(4,2)$};\node at (-0.5,0.75) {$\beta_{4,2}$};

\node at (1.65,0.7) {$\alpha_{1,2}$};
\node at (1.65,-0.7) {$\beta_{1,3}$};
\node at (-1.65,0.7) {$\beta_{3,3}$};
\node at (-1.65,-0.7) {$\alpha_{3,2}$};

\node at (0.7,1.65) {$\beta_{4,3}$};
\node at (0.7,-1.65) {$\alpha_{2,2}$};
\node at (-0.7,-1.65) {$\beta_{2,3}$};
\node at (-0.7,1.65) {$\alpha_{4,2}$};

\draw[arrows={-angle 90}, shorten >=15, shorten <=15]  (1,1) -- (2,0);
\draw[arrows={-angle 90}, shorten >=15, shorten <=15]  (1,-1) -- (0,-2);
\draw[arrows={-angle 90}, shorten >=15, shorten <=15]  (-1,-1) -- (-2,0);
\draw[arrows={-angle 90}, shorten >=15, shorten <=15]  (-1,1) -- (0,2);

\draw[arrows={-angle 90}, shorten >=15, shorten <=15]  (1,1) -- (1,0);
\draw[arrows={-angle 90}, shorten >=15, shorten <=15]  (1,-1) -- (0,-1);
\draw[arrows={-angle 90}, shorten >=15, shorten <=15]  (-1,-1) -- (-1,0);
\draw[arrows={-angle 90}, shorten >=15, shorten <=15]  (-1,1) -- (0,1);

\node at (0,2) {$(4,3)$};\node at (1,2.2) {$\alpha_{4,3}$};
\node at (2,0) {$(1,3)$};\node at (2.3,-1) {$\alpha_{1,3}$};
\node at (0,-2) {$(2,3)$};\node at (-1,-2.2) {$\alpha_{2,3}$};
\node at (-2,0) {$(3,3)$};\node at (-2.3,1) {$\alpha_{3,3}$};

\draw[arrows={-angle 90}, shorten >=15, shorten <=15]  (0,2) -- (2,2);
\draw[arrows={-angle 90}, shorten >=15, shorten <=15]  (2,0) -- (2,-2);
\draw[arrows={-angle 90}, shorten >=15, shorten <=15]  (0,-2) -- (-2,-2);
\draw[arrows={-angle 90}, shorten >=15, shorten <=15]  (-2,0) -- (-2,2);

\draw[arrows={-angle 90}, shorten >=15, shorten <=15]  (0,2) -- (1,1);
\draw[arrows={-angle 90}, shorten >=15, shorten <=15]  (2,0) -- (1,-1);
\draw[arrows={-angle 90}, shorten >=15, shorten <=15]  (0,-2) -- (-1,-1);
\draw[arrows={-angle 90}, shorten >=15, shorten <=15]  (-2,0) -- (-1,1);

\node at (2,2) {$(4,4)$};\node at (2.3,1) {$\beta_{4,4}$};
\node at (2,-2) {$(1,4)$};\node at (1,-2.2) {$\beta_{1,4}$};
\node at (-2,-2) {$(2,4)$};\node at (-2.3,-1) {$\beta_{2,4}$};
\node at (-2,2) {$(3,4)$};\node at (-1,2.2) {$\beta_{3,4}$};

\draw[arrows={-angle 90}, shorten >=15, shorten <=15]  (2,2) -- (2,0);
\draw[arrows={-angle 90}, shorten >=15, shorten <=15]  (2,-2) -- (0,-2);
\draw[arrows={-angle 90}, shorten >=15, shorten <=15]  (-2,-2) -- (-2,0);
\draw[arrows={-angle 90}, shorten >=15, shorten <=15]  (-2,2) -- (0,2);



\end{tikzpicture}
\end{center}
\end{ex}

We define a functor $ \Lambda : \mr{rep}_\C(\Delta_n,\mr{I}_N) \to \mr{rep}_\C(\hat{\Delta}_{n,N},\hat{\mr{I}}_{n,N})$ on objects by 
 \[\Lambda(M) := \hat{M} = \big( (\hat{M}^{(i,k)})_{i \in \Z_n, k \in [N] }, (\hat{M}_{\alpha_{i,k}}, \hat{M}_{\beta_{i,k+1}})_{i \in \Z_n, k \in [N-1]} \big)\]
 with 
 \begin{align*}
 \hat{M}^{(i,1)} &:= M^{(i)}   &\mr{for} \ k =1\\
\hat{M}^{(i,k)} &:= M_{i+k-2} \circ  M_{i+k-3} \circ  \dots \circ M_{i+1} \circ  M_{i} (M^{(i)})  &\mr{for} \ k \geq 2\\
\hat{M}_{\alpha_{i,k}} &:= M_{i+k-1}  &\mr{for} \ k \geq 1\\
\hat{M}_{\beta_{i,k}} &:= \iota : \hat{M}^{(i,k)}  \hookrightarrow \hat{M}^{(i+1,k-1)}  &\mr{for} \ k \geq 2
 \end{align*}
where the inclusions in the last row arise naturally from the definition of the vector spaces of  the representation $\hat{M}$.

\begin{ex}
Let $n=N=2$ and $M=U(1;2)\oplus U(2;2)$. The $\hat{\Delta}_{2,2}$-representation
$\Lambda(M)$ is
\[
\begin{tikzpicture}[scale=0.7]
\node at (0,2) {$\C$};
\node at (2,0) {$\C^2$};
\node at (0,-2) {$\C$};
\node at (-2,0) {$\C^2$};

\draw[arrows={-angle 90}, shorten >=9, shorten <=9]  (0,2) -- (2,0);
\draw[arrows={-angle 90}, shorten >=9, shorten <=9]  (2,0) -- (0,-2);
\draw[arrows={-angle 90}, shorten >=9, shorten <=9]  (0,-2) -- (-2,0);
\draw[arrows={-angle 90}, shorten >=9, shorten <=9]  (-2,0) -- (0,2);

\node at (1.5,1.5) {$\binom{0}{1}$};
\node at (1.5,-1.5) {$(1 0)$};
\node at (-1.5,-1.5) {$\binom{0}{1}$};
\node at (-1.5,1.5) {$(1 0)$};
\end{tikzpicture}
\]
\end{ex}

\begin{prop}\label{prop:Lambda-in-fully-faithful}
$ \Lambda : \mr{rep}_\C(\Delta_n,\mr{I}_N) \to \mr{rep}_\C(\hat{\Delta}_{n,N},\hat{\mr{I}}_{n,N})$ as defined above induces a bijection $\Lambda_{N,M} : \mr{Hom}_{\Delta_n}(N,M) \to \mr{Hom}_{\hat{\Delta}_{n,N}}(\hat{N},\hat{M})$ for all $N,M \in  \mr{rep}_\C(\Delta_n,\mr{I}_N) $ and hence is a fully faithful functor.
\end{prop}

\begin{proof}
By construction of $\Lambda$, the vector spaces constituting $\hat{M} \in \mr{rep}_\C(\hat{\Delta}_{n,N},\hat{\mr{I}}_{n,N})$ are subspaces of the corresponding vector spaces constituting $M$. Hence each morphism in $\mr{Hom}_{\Delta_n}(N,M)$ induces a morphism in $\mr{Hom}_{\hat{\Delta}_{n,N}}(\hat{N},\hat{M})$ whose components at the additional vertices are obtained by restriction. It is immediate to check that this induces the desired bijection  $\Lambda_{N,M}$,  $\Lambda_{N,M}(\mr{id}_M) = \mr{id}_{\hat{M}}$ and that $\Lambda_{N,M}(\phi) \circ \Lambda_{N,M} (\psi) = \Lambda_{N,M}(\phi \circ \psi)$ holds for all $\phi, \psi \in \mr{Hom}_{\Delta_n}(N,M)$ and all $N,M \in \mr{rep}_\C(\Delta_n,\mr{I}_N)$. 
\end{proof}

Now we want to describe the image of the indecomposable $U_i(\ell)$ under $\Lambda$. Let $A_{\infty \times N}$ be the infinite band quiver of height $N$, that is, the quiver with vertices $(i,k)$ for $i \in \Z$ and $k \in [N]$ and arrows $\alpha_{i,k} : (i,k) \to (i,k+1)$ and $\beta_{i,k} : (i,k) \to (i+1,k-1)$ whenever both vertices exist. Define a map of quivers $\phi : A_{\infty \times N} \to \hat{\Delta}_{n,N}$, induced by sending each index $i \in \Z$ to its equivalence class $i \in \Z_n$. This extends to a push-down functor $\Phi : \mr{rep}_\C(A_{\infty \times N}) \to \mr{rep}_\C(\hat{\Delta}_{n,N})$ with
\begin{align*}
\Big(\Phi (V)\Big)^{(i,k)} &= \bigoplus_{r \in \Z} V^{(i+rn,k)}\\ 
\Big(\Phi (V) \Big)_{\alpha_{i,k}} &= \bigoplus_{r \in \Z} V_{\alpha_{i+rn,k}}\\
\Big(\Phi(V) \Big)_{\beta_{i,k}} &= \bigoplus_{r \in \Z} V_{\beta_{i+rn,k}}
\end{align*}
for all $V \in \mr{rep}_\C(A_{\infty \times N})$. Consider the $A_{\infty \times N}$-representation $V({i;\ell})$ with vector spaces $V({i;\ell})^{(j,k)} = \C$ for $(j,k) \in (A_{\infty \times N})_0$ with $i \leq j \leq i+ \ell -1$ and $1 \leq k \leq i+ \ell - j$ and zero otherwise. The maps along the arrows of $A_{\infty \times N}$ are identities if both the source and target space are one-dimensional and zero otherwise. Using the explicit definitions of the functors $\Lambda$ and $\Phi$, and the explicit descriptions of the representations $U_i(\ell)$ and $V(i;\ell)$, we can now directly verify that $$\Lambda( U_i(\ell) ) = \Phi ( V(i;\ell) ).$$

Analogously, we define  $A_{\infty \times N}$-representations $V({i,k;\ell})$ consisting of vector spaces $V({i,k;\ell})^{(j,r)} = \C$ for $(j,r) \in (A_{\infty \times N})_0$ with $j\geq i$ and $i+k \leq j+r \leq i+ k+ \ell$ and $W({i,k;\ell})$ with vector spaces $W({i,k;\ell})^{(j,r)} = \C$ for $(j,r) \in (A_{\infty \times N})_0$ with $j\leq i$ and $i-\ell+1 \leq j+r \leq i+ k$.

\begin{ex}
For $N=4$ the quiver $A_{\infty \times N}$ is 
\[
\begin{tikzpicture}[scale=1.0]
\node at (-5,0) {$\dots$};
\node at (-4,0) {\tiny$(-4,1)$};
\node at (-3,0) {\tiny$(-3,1)$};
\node at (-2,0) {\tiny$(-2,1)$};
\node at (-1,0) {\tiny$(-1,1)$};
\node at (0,0) {\tiny$(0,1)$};
\node at (1,0) {\tiny$(1,1)$};
\node at (2,0) {\tiny$(2,1)$};
\node at (3,0) {\tiny$(3,1)$};
\node at (4,0) {\tiny$(4,1)$};
\node at (5,0) {$\dots$};

\draw[arrows={-angle 90}, shorten >=9, shorten <=9]  (-4,0) -- (-4+0.5,1);
\draw[arrows={-angle 90}, shorten >=9, shorten <=9]  (-4+0.5,1) -- (-3,0);

\draw[arrows={-angle 90}, shorten >=9, shorten <=9]  (-3,0) -- (-3+0.5,1);
\draw[arrows={-angle 90}, shorten >=9, shorten <=9]  (-3+0.5,1) -- (-3+1,0);

\draw[arrows={-angle 90}, shorten >=9, shorten <=9]  (-2,0) -- (-2+0.5,1);
\draw[arrows={-angle 90}, shorten >=9, shorten <=9]  (-2+0.5,1) -- (-2+1,0);

\draw[arrows={-angle 90}, shorten >=9, shorten <=9]  (-1,0) -- (-1+0.5,1);
\draw[arrows={-angle 90}, shorten >=9, shorten <=9]  (-1+0.5,1) -- (-1+1,0);

\draw[arrows={-angle 90}, shorten >=9, shorten <=9]  (-0,0) -- (-0+0.5,1);
\draw[arrows={-angle 90}, shorten >=9, shorten <=9]  (-0+0.5,1) -- (-0+1,0);

\draw[arrows={-angle 90}, shorten >=9, shorten <=9]  (1,0) -- (1+0.5,1);
\draw[arrows={-angle 90}, shorten >=9, shorten <=9]  (1+0.5,1) -- (1+1,0);

\draw[arrows={-angle 90}, shorten >=9, shorten <=9]  (2,0) -- (2+0.5,1);
\draw[arrows={-angle 90}, shorten >=9, shorten <=9]  (2+0.5,1) -- (2+1,0);

\draw[arrows={-angle 90}, shorten >=9, shorten <=9]  (3,0) -- (3+0.5,1);
\draw[arrows={-angle 90}, shorten >=9, shorten <=9]  (3+0.5,1) -- (3+1,0);

\draw[arrows={-angle 90}, shorten >=9, shorten <=9]  (4,0) -- (4+0.5,1);

\node at (-5+0.5,1) {$\dots$};
\node at (-4+0.5,1) {\tiny$(-4,2)$};
\node at (-3+0.5,1) {\tiny$(-3,2)$};
\node at (-2+0.5,1) {\tiny$(-2,2)$};
\node at (-1+0.5,1) {\tiny$(-1,2)$};
\node at (0+0.5,1) {\tiny$(0,2)$};
\node at (1+0.5,1) {\tiny$(1,2)$};
\node at (2+0.5,1) {\tiny$(2,2)$};
\node at (3+0.5,1) {\tiny$(3,2)$};
\node at (4+0.5,1) {\tiny$(4,2)$};
\node at (5+0.5,1) {$\dots$};

\draw[arrows={-angle 90}, shorten >=9, shorten <=9]  (-4,2) -- (-4+0.5,1);
\draw[arrows={-angle 90}, shorten >=9, shorten <=9]  (-4+0.5,1) -- (-3,2);

\draw[arrows={-angle 90}, shorten >=9, shorten <=9]  (-3,2) -- (-3+0.5,1);
\draw[arrows={-angle 90}, shorten >=9, shorten <=9]  (-3+0.5,1) -- (-3+1,2);

\draw[arrows={-angle 90}, shorten >=9, shorten <=9]  (-2,2) -- (-2+0.5,1);
\draw[arrows={-angle 90}, shorten >=9, shorten <=9]  (-2+0.5,1) -- (-2+1,2);

\draw[arrows={-angle 90}, shorten >=9, shorten <=9]  (-1,2) -- (-1+0.5,1);
\draw[arrows={-angle 90}, shorten >=9, shorten <=9]  (-1+0.5,1) -- (-1+1,2);

\draw[arrows={-angle 90}, shorten >=9, shorten <=9]  (-0,2) -- (-0+0.5,1);
\draw[arrows={-angle 90}, shorten >=9, shorten <=9]  (-0+0.5,1) -- (-0+1,2);

\draw[arrows={-angle 90}, shorten >=9, shorten <=9]  (1,2) -- (1+0.5,1);
\draw[arrows={-angle 90}, shorten >=9, shorten <=9]  (1+0.5,1) -- (1+1,2);

\draw[arrows={-angle 90}, shorten >=9, shorten <=9]  (2,2) -- (2+0.5,1);
\draw[arrows={-angle 90}, shorten >=9, shorten <=9]  (2+0.5,1) -- (2+1,2);

\draw[arrows={-angle 90}, shorten >=9, shorten <=9]  (3,2) -- (3+0.5,1);
\draw[arrows={-angle 90}, shorten >=9, shorten <=9]  (3+0.5,1) -- (3+1,2);

\draw[arrows={-angle 90}, shorten >=9, shorten <=9]  (4,2) -- (4+0.5,1);

\node at (-5,2) {$\dots$};
\node at (-4,2) {\tiny$(-5,3)$};
\node at (-3,2) {\tiny$(-4,3)$};
\node at (-2,2) {\tiny$(-3,3)$};
\node at (-1,2) {\tiny$(-2,3)$};
\node at (0,2) {\tiny$(-1,3)$};
\node at (1,2) {\tiny$(0,3)$};
\node at (2,2) {\tiny$(1,3)$};
\node at (3,2) {\tiny$(2,3)$};
\node at (4,2) {\tiny$(3,3)$};
\node at (5,2) {$\dots$};

\draw[arrows={-angle 90}, shorten >=9, shorten <=9]  (-4,2) -- (-4+0.5,3);
\draw[arrows={-angle 90}, shorten >=9, shorten <=9]  (-4+0.5,3) -- (-3,2);

\draw[arrows={-angle 90}, shorten >=9, shorten <=9]  (-3,2) -- (-3+0.5,3);
\draw[arrows={-angle 90}, shorten >=9, shorten <=9]  (-3+0.5,3) -- (-3+1,2);

\draw[arrows={-angle 90}, shorten >=9, shorten <=9]  (-2,2) -- (-2+0.5,3);
\draw[arrows={-angle 90}, shorten >=9, shorten <=9]  (-2+0.5,3) -- (-2+1,2);

\draw[arrows={-angle 90}, shorten >=9, shorten <=9]  (-1,2) -- (-1+0.5,3);
\draw[arrows={-angle 90}, shorten >=9, shorten <=9]  (-1+0.5,3) -- (-1+1,2);

\draw[arrows={-angle 90}, shorten >=9, shorten <=9]  (-0,2) -- (-0+0.5,3);
\draw[arrows={-angle 90}, shorten >=9, shorten <=9]  (-0+0.5,3) -- (-0+1,2);

\draw[arrows={-angle 90}, shorten >=9, shorten <=9]  (1,2) -- (1+0.5,3);
\draw[arrows={-angle 90}, shorten >=9, shorten <=9]  (1+0.5,3) -- (1+1,2);

\draw[arrows={-angle 90}, shorten >=9, shorten <=9]  (2,2) -- (2+0.5,3);
\draw[arrows={-angle 90}, shorten >=9, shorten <=9]  (2+0.5,3) -- (2+1,2);

\draw[arrows={-angle 90}, shorten >=9, shorten <=9]  (3,2) -- (3+0.5,3);
\draw[arrows={-angle 90}, shorten >=9, shorten <=9]  (3+0.5,3) -- (3+1,2);

\draw[arrows={-angle 90}, shorten >=9, shorten <=9]  (4,2) -- (4+0.5,3);

\node at (-5+0.5,3) {$\dots$};
\node at (-4+0.5,3) {\tiny$(-5,4)$};
\node at (-3+0.5,3) {\tiny$(-4,4)$};
\node at (-2+0.5,3) {\tiny$(-3,4)$};
\node at (-1+0.5,3) {\tiny$(-2,4)$};
\node at (0+0.5,3) {\tiny$(-1,4)$};
\node at (1+0.5,3) {\tiny$(0,4)$};
\node at (2+0.5,3) {\tiny$(1,4)$};
\node at (3+0.5,3) {\tiny$(2,4)$};
\node at (4+0.5,3) {\tiny$(3,4)$};
\node at (5+0.5,3) {$\dots$};

\end{tikzpicture}
\]
Its representations $V(i,3)$ are of the form
\[
\begin{tikzpicture}[scale=1.0]
\node at (-5,0) {$\dots$};
\node at (-4,0) {$0$};
\node at (-3,0) {$0$};
\node at (-2,0) {$0$};
\node at (-1,0) {$0$};
\node at (0,0) {$\C$};
\node at (1,0) {$\C$};
\node at (2,0) {$\C$};
\node at (3,0) {$0$};
\node at (4,0) {$0$};
\node at (5,0) {$\dots$};

\draw[arrows={-angle 90}, shorten >=9, shorten <=9]  (-4,0) -- (-4+0.5,1);
\draw[arrows={-angle 90}, shorten >=9, shorten <=9]  (-4+0.5,1) -- (-3,0);

\draw[arrows={-angle 90}, shorten >=9, shorten <=9]  (-3,0) -- (-3+0.5,1);
\draw[arrows={-angle 90}, shorten >=9, shorten <=9]  (-3+0.5,1) -- (-3+1,0);

\draw[arrows={-angle 90}, shorten >=9, shorten <=9]  (-2,0) -- (-2+0.5,1);
\draw[arrows={-angle 90}, shorten >=9, shorten <=9]  (-2+0.5,1) -- (-2+1,0);

\draw[arrows={-angle 90}, shorten >=9, shorten <=9]  (-1,0) -- (-1+0.5,1);
\draw[arrows={-angle 90}, shorten >=9, shorten <=9]  (-1+0.5,1) -- (-1+1,0);

\draw[arrows={-angle 90}, shorten >=9, shorten <=9]  (-0,0) -- (-0+0.5,1);
\draw[arrows={-angle 90}, shorten >=9, shorten <=9]  (-0+0.5,1) -- (-0+1,0);

\draw[arrows={-angle 90}, shorten >=9, shorten <=9]  (1,0) -- (1+0.5,1);
\draw[arrows={-angle 90}, shorten >=9, shorten <=9]  (1+0.5,1) -- (1+1,0);

\draw[arrows={-angle 90}, shorten >=9, shorten <=9]  (2,0) -- (2+0.5,1);
\draw[arrows={-angle 90}, shorten >=9, shorten <=9]  (2+0.5,1) -- (2+1,0);

\draw[arrows={-angle 90}, shorten >=9, shorten <=9]  (3,0) -- (3+0.5,1);
\draw[arrows={-angle 90}, shorten >=9, shorten <=9]  (3+0.5,1) -- (3+1,0);

\draw[arrows={-angle 90}, shorten >=9, shorten <=9]  (4,0) -- (4+0.5,1);

\node at (-5+0.5,1) {$\dots$};
\node at (-4+0.5,1) {$0$};
\node at (-3+0.5,1) {$0$};
\node at (-2+0.5,1) {$0$};
\node at (-1+0.5,1) {$0$};
\node at (0+0.5,1) {$\C$};
\node at (1+0.5,1) {$\C$};
\node at (2+0.5,1) {$0$};
\node at (3+0.5,1) {$0$};
\node at (4+0.5,1) {$0$};
\node at (5+0.5,1) {$\dots$};

\draw[arrows={-angle 90}, shorten >=9, shorten <=9]  (-4,2) -- (-4+0.5,1);
\draw[arrows={-angle 90}, shorten >=9, shorten <=9]  (-4+0.5,1) -- (-3,2);

\draw[arrows={-angle 90}, shorten >=9, shorten <=9]  (-3,2) -- (-3+0.5,1);
\draw[arrows={-angle 90}, shorten >=9, shorten <=9]  (-3+0.5,1) -- (-3+1,2);

\draw[arrows={-angle 90}, shorten >=9, shorten <=9]  (-2,2) -- (-2+0.5,1);
\draw[arrows={-angle 90}, shorten >=9, shorten <=9]  (-2+0.5,1) -- (-2+1,2);

\draw[arrows={-angle 90}, shorten >=9, shorten <=9]  (-1,2) -- (-1+0.5,1);
\draw[arrows={-angle 90}, shorten >=9, shorten <=9]  (-1+0.5,1) -- (-1+1,2);

\draw[arrows={-angle 90}, shorten >=9, shorten <=9]  (-0,2) -- (-0+0.5,1);
\draw[arrows={-angle 90}, shorten >=9, shorten <=9]  (-0+0.5,1) -- (-0+1,2);

\draw[arrows={-angle 90}, shorten >=9, shorten <=9]  (1,2) -- (1+0.5,1);
\draw[arrows={-angle 90}, shorten >=9, shorten <=9]  (1+0.5,1) -- (1+1,2);

\draw[arrows={-angle 90}, shorten >=9, shorten <=9]  (2,2) -- (2+0.5,1);
\draw[arrows={-angle 90}, shorten >=9, shorten <=9]  (2+0.5,1) -- (2+1,2);

\draw[arrows={-angle 90}, shorten >=9, shorten <=9]  (3,2) -- (3+0.5,1);
\draw[arrows={-angle 90}, shorten >=9, shorten <=9]  (3+0.5,1) -- (3+1,2);

\draw[arrows={-angle 90}, shorten >=9, shorten <=9]  (4,2) -- (4+0.5,1);

\node at (-5,2) {$\dots$};
\node at (-4,2) {$0$};
\node at (-3,2) {$0$};
\node at (-2,2) {$0$};
\node at (-1,2) {$0$};
\node at (0,2) {$0$};
\node at (1,2) {$\C$};
\node at (2,2) {$0$};
\node at (3,2) {$0$};
\node at (4,2) {$0$};
\node at (5,2) {$\dots$};

\draw[arrows={-angle 90}, shorten >=9, shorten <=9]  (-4,2) -- (-4+0.5,3);
\draw[arrows={-angle 90}, shorten >=9, shorten <=9]  (-4+0.5,3) -- (-3,2);

\draw[arrows={-angle 90}, shorten >=9, shorten <=9]  (-3,2) -- (-3+0.5,3);
\draw[arrows={-angle 90}, shorten >=9, shorten <=9]  (-3+0.5,3) -- (-3+1,2);

\draw[arrows={-angle 90}, shorten >=9, shorten <=9]  (-2,2) -- (-2+0.5,3);
\draw[arrows={-angle 90}, shorten >=9, shorten <=9]  (-2+0.5,3) -- (-2+1,2);

\draw[arrows={-angle 90}, shorten >=9, shorten <=9]  (-1,2) -- (-1+0.5,3);
\draw[arrows={-angle 90}, shorten >=9, shorten <=9]  (-1+0.5,3) -- (-1+1,2);

\draw[arrows={-angle 90}, shorten >=9, shorten <=9]  (-0,2) -- (-0+0.5,3);
\draw[arrows={-angle 90}, shorten >=9, shorten <=9]  (-0+0.5,3) -- (-0+1,2);

\draw[arrows={-angle 90}, shorten >=9, shorten <=9]  (1,2) -- (1+0.5,3);
\draw[arrows={-angle 90}, shorten >=9, shorten <=9]  (1+0.5,3) -- (1+1,2);

\draw[arrows={-angle 90}, shorten >=9, shorten <=9]  (2,2) -- (2+0.5,3);
\draw[arrows={-angle 90}, shorten >=9, shorten <=9]  (2+0.5,3) -- (2+1,2);

\draw[arrows={-angle 90}, shorten >=9, shorten <=9]  (3,2) -- (3+0.5,3);
\draw[arrows={-angle 90}, shorten >=9, shorten <=9]  (3+0.5,3) -- (3+1,2);

\draw[arrows={-angle 90}, shorten >=9, shorten <=9]  (4,2) -- (4+0.5,3);

\node at (-5+0.5,3) {$\dots$};
\node at (-4+0.5,3) {$0$};
\node at (-3+0.5,3) {$0$};
\node at (-2+0.5,3) {$0$};
\node at (-1+0.5,3) {$0$};
\node at (0+0.5,3) {$0$};
\node at (1+0.5,3) {$0$};
\node at (2+0.5,3) {$0$};
\node at (3+0.5,3) {$0$};
\node at (4+0.5,3) {$0$};
\node at (5+0.5,3) {$\dots$};

\end{tikzpicture}
\]
From now on we erase all zeros and arrows connected to zeros from the picture. Hence we obtain
\[
\begin{tikzpicture}[scale=1.0]
\node at (-1,1) {$V(i,2,2)=$};

\node at (1,0) {$\C$};
\node at (2,0) {$\C$};

\draw[arrows={-angle 90}, shorten >=9, shorten <=9]  (-0+0.5,1) -- (-0+1,0);

\draw[arrows={-angle 90}, shorten >=9, shorten <=9]  (1,0) -- (1+0.5,1);
\draw[arrows={-angle 90}, shorten >=9, shorten <=9]  (1+0.5,1) -- (1+1,0);

\node at (0+0.5,1) {$\C$};
\node at (1+0.5,1) {$\C$};

\draw[arrows={-angle 90}, shorten >=9, shorten <=9]  (-0+0.5,1) -- (-0+1,2);

\draw[arrows={-angle 90}, shorten >=9, shorten <=9]  (1,2) -- (1+0.5,1);

\node at (1,2) {$\C$};

\end{tikzpicture}
\quad \quad
\begin{tikzpicture}[scale=1.0]
\node at (-3.5,1.5) {$\mr{and}$};
\node at (-1,1.5) {$W(i,3,2)=$};
\node at (0,0) {$\C$};
\node at (1,0) {$\C$};

\draw[arrows={-angle 90}, shorten >=9, shorten <=9]  (-0,0) -- (-0+0.5,1);
\draw[arrows={-angle 90}, shorten >=9, shorten <=9]  (-0+0.5,1) -- (-0+1,0);
\draw[arrows={-angle 90}, shorten >=9, shorten <=9]  (1,0) -- (1+0.5,1);

\node at (0+0.5,1) {$\C$};
\node at (1+0.5,1) {$\C$};

\draw[arrows={-angle 90}, shorten >=9, shorten <=9]  (-0+0.5,1) -- (-0+1,2);

\draw[arrows={-angle 90}, shorten >=9, shorten <=9]  (1,2) -- (1+0.5,1);
\draw[arrows={-angle 90}, shorten >=9, shorten <=9]  (1+0.5,1) -- (1+1,2);

\node at (1,2) {$\C$};
\node at (2,2) {$\C$};

\draw[arrows={-angle 90}, shorten >=9, shorten <=9]  (1,2) -- (1+0.5,3);
\draw[arrows={-angle 90}, shorten >=9, shorten <=9]  (1+0.5,3) -- (1+1,2);

\node at (1+0.5,3) {$\C$};

\end{tikzpicture}
\]
\end{ex}

\subsection{Homological Properties of the Category of Cylinder Representations}

In this section, we follow closely the approach of \cite{CFR2013} to establish certain favourable homological properties of the image of the functor $\Lambda$.

\begin{prop}\label{prop:proj-and-inj-objects}
The simple, projective and injective objects in $ \mr{rep}_\C(\hat{\Delta}_{n,N},\hat{\mr{I}}_{n,N})$ are given as
\begin{align*}
 S_{i,k} := \Phi\big(V(i,k;0)\big), \quad  P_{i,k} := \Phi\big(V(i,k;N-k)\big),  \quad I_{i,k} := \Phi\big(W(i,k;N-k)\big),
 \end{align*}
respectively, for all  $(i,k) \in (\hat{\Delta}_{n,N})_0$.
\end{prop}

\begin{proof}
For the simple objects this is immediate. The parametrization of the projective and injective representations is a direct computation using the formula based on paths in the quiver $\hat{\Delta}_{n,N}$ (see  \cite[Definition~5.3]{Schiffler2014}) and their relations from $\hat{\mr{I}}_{n,N}$ as described in the beginning of this section.
\end{proof}

\begin{trm}\label{trm:global-dim-two}
The category $\mr{rep}_\C(\hat{\Delta}_{n,N},\hat{\mr{I}}_{n,N})$ has global dimension at most two.
\end{trm}

\begin{proof}
It suffices to construct projective resolutions of length at most two for all simple representations in $\mr{rep}_\C(\hat{\Delta}_{n,N},\hat{\mr{I}}_{n,N})$. 
These representations are denoted by $S_{i,k}$ and consist of a single copy of $\C$ at vertex $(i,k)$ and all other vector spaces and the maps are zero. The projective resolutions of  $S_{i,1}$ are of the form
\[ 0 \to P_{i,2} \to P_{i,1} \to S_{i,1} \to 0\]
and for $S_{i,k}$ with $k \geq 2$ this generalizes to 
\[ 0 \to P_{i+1,k} \to P_{i+1,k-1} \oplus P_{i,k+1} \to P_{i,k} \to S_{i,k} \to 0. \]
\end{proof}
\begin{ex}
For $N=4$ and $S_{i,3}$ we obtain the following projective resolution:
\[
\begin{tikzpicture}[scale=1.0]
\node at (-5,0) {$\C$};
\node at (-4,0) {$\C$};
\node at (-5-0.5,1) {$\C$};
\node at (-4-0.5,1) {$\C$};
\node at (-6,2) {$\C$};
\node at (-5,2) {$\C$};
\node at (-5-0.5,3) {$\C$};
\draw[arrows={-angle 90}, shorten >=9, shorten <=9]  (-6+0.5,1) -- (-6+1,0);
\draw[arrows={-angle 90}, shorten >=9, shorten <=9]  (-5,0) -- (-5+0.5,1);
\draw[arrows={-angle 90}, shorten >=9, shorten <=9]  (-5+0.5,1) -- (-4,0);
\draw[arrows={-angle 90}, shorten >=9, shorten <=9]  (-6,2) -- (-6+0.5,1);
\draw[arrows={-angle 90}, shorten >=9, shorten <=9]  (-5-0.5,1) -- (-5,2);
\draw[arrows={-angle 90}, shorten >=9, shorten <=9]  (-5,2) -- (-4-0.5,1);
\draw[arrows={-angle 90}, shorten >=9, shorten <=9]  (-6,2) -- (-5-0.5,3);
\draw[arrows={-angle 90}, shorten >=9, shorten <=9]  (-5-0.5,3) -- (-5,2);

\node at (-2,0) {$\C$};
\node at (-1,0) {$\C^2$};
\node at (0,0) {$\C$};
\node at (2,0) {$\C$};
\node at (3,0) {$\C$};

\draw[arrows={-angle 90}, shorten >=9, shorten <=9]  (-3+0.5,1) -- (-3+1,0);

\draw[arrows={-angle 90}, shorten >=9, shorten <=9]  (-2,0) -- (-2+0.5,1);
\draw[arrows={-angle 90}, shorten >=9, shorten <=9]  (-2+0.5,1) -- (-2+1,0);

\draw[arrows={-angle 90}, shorten >=9, shorten <=9]  (-1,0) -- (-1+0.5,1);
\draw[arrows={-angle 90}, shorten >=9, shorten <=9]  (-1+0.5,1) -- (-1+1,0);


\draw[arrows={-angle 90}, shorten >=9, shorten <=9]  (1+0.5,1) -- (1+1,0);

\draw[arrows={-angle 90}, shorten >=9, shorten <=9]  (2,0) -- (2+0.5,1);
\draw[arrows={-angle 90}, shorten >=9, shorten <=9]  (2+0.5,1) -- (2+1,0);

\node at (-3+0.5,1) {$\C$};
\node at (-2+0.5,1) {$\C^2$};
\node at (-1+0.5,1) {$\C$};
\node at (1+0.5,1) {$\C$};
\node at (2+0.5,1) {$\C$};

\draw[arrows={-angle 90}, shorten >=9, shorten <=9]  (-3+0.5,1) -- (-3+1,2);

\draw[arrows={-angle 90}, shorten >=9, shorten <=9]  (-2,2) -- (-2+0.5,1);
\draw[arrows={-angle 90}, shorten >=9, shorten <=9]  (-2+0.5,1) -- (-2+1,2);

\draw[arrows={-angle 90}, shorten >=9, shorten <=9]  (-1,2) -- (-1+0.5,1);


\draw[arrows={-angle 90}, shorten >=9, shorten <=9]  (1,2) -- (1+0.5,1);
\draw[arrows={-angle 90}, shorten >=9, shorten <=9]  (1+0.5,1) -- (1+1,2);

\draw[arrows={-angle 90}, shorten >=9, shorten <=9]  (2,2) -- (2+0.5,1);

\node at (-2,2) {$\C^2$};
\node at (-1,2) {$\C$};
\node at (1,2) {$\C$};
\node at (2,2) {$\C$};

\node at (4,2) {$\C$};

\draw[arrows={-angle 90}, shorten >=9, shorten <=9]  (-3+0.5,3) -- (-3+1,2);

\draw[arrows={-angle 90}, shorten >=9, shorten <=9]  (-2,2) -- (-2+0.5,3);
\draw[arrows={-angle 90}, shorten >=9, shorten <=9]  (-2+0.5,3) -- (-2+1,2);



\draw[arrows={-angle 90}, shorten >=9, shorten <=9]  (1,2) -- (1+0.5,3);
\draw[arrows={-angle 90}, shorten >=9, shorten <=9]  (1+0.5,3) -- (1+1,2);

\node at (-3+0.5,3) {$\C$};
\node at (-2+0.5,3) {$\C$};
\node at (1+0.5,3) {$\C$};

\draw[arrows={-angle 90}, shorten >=5, shorten <=5]  (-0.5,2) -- (0.5,2);

\draw[arrows={-angle 90}, shorten >=5, shorten <=5,->>]  (2.5,2) -- (3.5,2);

\draw[arrows={-angle 90}, shorten >=5, shorten <=5,right hook ->]  (-4,2) -- (-3,2);

\node at (4,4) {$S_{i,3}$};
\node at (1.5,4) {$P_{i,3}$};
\node at (-2,4) {$P_{1+1,2}\oplus P_{i,4}$};
\node at (-5.5,4) {$P_{i+1,3}$};
\end{tikzpicture}
\]
\end{ex}
\begin{lma}\label{lma:proj-inj-dim-one}
For $M \in \mr{rep}_\C(\Delta_n,\mr{I}_N)$ the injective and projective dimension of $\hat{M}$ is at most one and $\mr{Ext}^1_{\hat{\Delta}_{n,N},\hat{\mr{I}}_{n,N}}(\hat{M},\hat{M}) = 0$.
\end{lma}

\begin{proof}
It suffices to compute the projective and injective dimension of the image of all indecomposable representations $U_i(\ell) \in \mr{rep}_\C(\Delta_n,\mr{I}_N)$, by exhibiting projective resp.~injective resolutions, namely: %
\[ 0 \to P_{i,\ell+1} \to P_{i,1} \to \hat{U}_i(\ell) \to 0, \]
\[ 0 \to \hat{U}(j;\ell) \to I_{j,1} \to I_{j-\ell,\ell+1} \to 0 \]
where $j := i + \ell -1$ and hence $U(j;\ell) = U_i(\ell)$. 

It remains to prove vanishing of all $ \mr{Ext}^1_{\hat{\Delta}_{n,N},\hat{\mr{I}}_{n,N}} \big(   \hat{U}_i(\ell)  ,  \hat{U}_j(\ell') \big)$. We apply the functor ${\rm Hom_{\hat{\Delta}_{n,N},\hat{\mr{I}}_{n,N}}}(\_,\hat{U}_j(\ell'))$ to the above projective resolution of $\hat{U}_i(\ell)$, simplify the terms involving projectives, and obtain the exact sequence
$$0\rightarrow{\rm Hom}(\hat{U}_i(\ell),\hat{U}_j(\ell'))\rightarrow \hat{U}_j(\ell')^{(i,1)}\stackrel{\alpha}{\rightarrow}\hat{U}_j(\ell')^{(i,\ell+1)}\rightarrow{\rm Ext}^1(\hat{U}_i(\ell),\hat{U}_j(\ell'))\rightarrow 0.$$
By definition of $\Lambda$, the map $\alpha$ is the canonical surjection
$$U_j(\ell')^{(i)}\rightarrow(U_j(\ell')_{i+\ell-1}\circ\ldots\circ U_j(\ell')_i)(U_j(\ell')^{(i)},$$
proving the desired ${\rm Ext^1}$-vanishing.
\end{proof}
\begin{ex}
For $N=4$ we obtain the following projective and injective resolutions of ${U}_i(3)=U(i-2;3)$:
\[
\begin{tikzpicture}[scale=1.0]
\foreach \dist in {-5}{
\node at (0-0.5+\dist,4) {$P_{i,4}$};

\node at (1+\dist,0) {$\C$};
\draw[arrows={-angle 90}, shorten >=9, shorten <=9]  (0+0.5+\dist,1) -- (1+\dist,0);
\node at (0+0.5+\dist,1) {$\C$};
\draw[arrows={-angle 90}, shorten >=9, shorten <=9]  (-1+1+\dist,2) -- (0+0.5+\dist,1);
\node at (0+\dist,2) {$\C$};
\draw[arrows={-angle 90}, shorten >=9, shorten <=9]  (-1+0.5+\dist,3) -- (-0+\dist,2);
\node at (-1+0.5+\dist,3) {$\C$};

}

\foreach \dist in {-3}{
\draw[arrows={-angle 90}, shorten >=5, shorten <=5,right hook ->]  (-0.5+\dist,2) -- (0.5+\dist,2);
}

\foreach \dist in {0}{
\node at (0-0.5+\dist,4) {$P_{i,1}$};

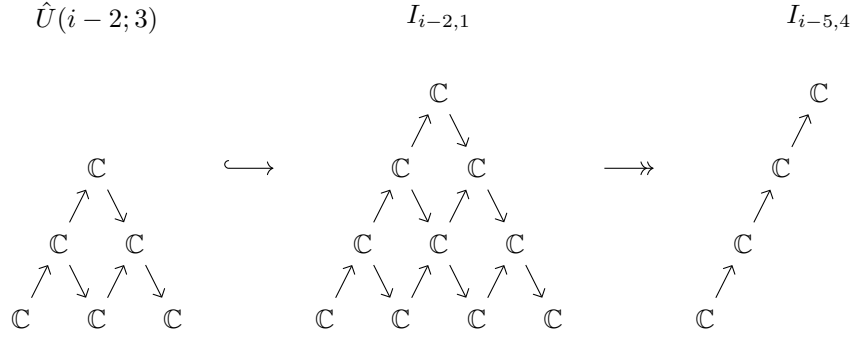
\node at (-2+\dist,0) {$\C$};
\node at (-1+\dist,0) {$\C$};
\node at (0+\dist,0) {$\C$};
\node at (1+\dist,0) {$\C$};

\draw[arrows={-angle 90}, shorten >=9, shorten <=9]  (-2+\dist,0) -- (-2+0.5+\dist,1);
\draw[arrows={-angle 90}, shorten >=9, shorten <=9]  (-2+0.5+\dist,1) -- (-1+\dist,0);
\draw[arrows={-angle 90}, shorten >=9, shorten <=9]  (-1+\dist,0) -- (-1+0.5+\dist,1);
\draw[arrows={-angle 90}, shorten >=9, shorten <=9]  (-1+0.5+\dist,1) -- (0+\dist,0);
\draw[arrows={-angle 90}, shorten >=9, shorten <=9]  (0+\dist,0) -- (0+0.5+\dist,1);
\draw[arrows={-angle 90}, shorten >=9, shorten <=9]  (0+0.5+\dist,1) -- (1+\dist,0);

\node at (-2+0.5+\dist,1) {$\C$};
\node at (-1+0.5+\dist,1) {$\C$};
\node at (0+0.5+\dist,1) {$\C$};

\draw[arrows={-angle 90}, shorten >=9, shorten <=9]  (-2+0.5+\dist,1) -- (-2+1+\dist,2);
\draw[arrows={-angle 90}, shorten >=9, shorten <=9]  (-2+1+\dist,2) -- (-1+0.5+\dist,1);
\draw[arrows={-angle 90}, shorten >=9, shorten <=9]  (-1+0.5+\dist,1) -- (-1+1+\dist,2);
\draw[arrows={-angle 90}, shorten >=9, shorten <=9]  (-1+1+\dist,2) -- (0+0.5+\dist,1);

\node at (-1+\dist,2) {$\C$};
\node at (0+\dist,2) {$\C$};

\draw[arrows={-angle 90}, shorten >=9, shorten <=9]  (-1+\dist,2) -- (-1+0.5+\dist,3);
\draw[arrows={-angle 90}, shorten >=9, shorten <=9]  (-1+0.5+\dist,3) -- (-0+\dist,2);

\node at (-1+0.5+\dist,3) {$\C$};

}

\foreach \dist in {2}{
\draw[arrows={-angle 90}, shorten >=5, shorten <=5,->>]  (-0.5+\dist,2) -- (0.5+\dist,2);
}

\foreach \dist in {4}{
\node at (0+\dist,4) {$\hat{U}_{i}(3)$};

\node at (-1+\dist,0) {$\C$};
\node at (0+\dist,0) {$\C$};
\node at (1+\dist,0) {$\C$};

\draw[arrows={-angle 90}, shorten >=9, shorten <=9]  (-1+\dist,0) -- (-1+0.5+\dist,1);
\draw[arrows={-angle 90}, shorten >=9, shorten <=9]  (-1+0.5+\dist,1) -- (0+\dist,0);
\draw[arrows={-angle 90}, shorten >=9, shorten <=9]  (0+\dist,0) -- (0+0.5+\dist,1);
\draw[arrows={-angle 90}, shorten >=9, shorten <=9]  (0+0.5+\dist,1) -- (1+\dist,0);

\node at (-1+0.5+\dist,1) {$\C$};
\node at (0+0.5+\dist,1) {$\C$};

\draw[arrows={-angle 90}, shorten >=9, shorten <=9]  (-1+0.5+\dist,1) -- (-1+1+\dist,2);
\draw[arrows={-angle 90}, shorten >=9, shorten <=9]  (-1+1+\dist,2) -- (0+0.5+\dist,1);

\node at (0+\dist,2) {$\C$};
}

\end{tikzpicture}
\]
\[
\begin{tikzpicture}[scale=1.0]
\foreach \dist in {-5}{
\node at (0+\dist,4) {$\hat{U}(i-2;3)$};

\node at (-1+\dist,0) {$\C$};
\node at (0+\dist,0) {$\C$};
\node at (1+\dist,0) {$\C$};

\draw[arrows={-angle 90}, shorten >=9, shorten <=9]  (-1+\dist,0) -- (-1+0.5+\dist,1);
\draw[arrows={-angle 90}, shorten >=9, shorten <=9]  (-1+0.5+\dist,1) -- (0+\dist,0);
\draw[arrows={-angle 90}, shorten >=9, shorten <=9]  (0+\dist,0) -- (0+0.5+\dist,1);
\draw[arrows={-angle 90}, shorten >=9, shorten <=9]  (0+0.5+\dist,1) -- (1+\dist,0);

\node at (-1+0.5+\dist,1) {$\C$};
\node at (0+0.5+\dist,1) {$\C$};

\draw[arrows={-angle 90}, shorten >=9, shorten <=9]  (-1+0.5+\dist,1) -- (-1+1+\dist,2);
\draw[arrows={-angle 90}, shorten >=9, shorten <=9]  (-1+1+\dist,2) -- (0+0.5+\dist,1);

\node at (0+\dist,2) {$\C$};
}

\foreach \dist in {-3}{
\draw[arrows={-angle 90}, shorten >=5, shorten <=5,right hook ->]  (-0.5+\dist,2) -- (0.5+\dist,2);
}

\foreach \dist in {0}{
\node at (0-0.5+\dist,4) {$I_{i-2,1}$};

\node at (-2+\dist,0) {$\C$};
\node at (-1+\dist,0) {$\C$};
\node at (0+\dist,0) {$\C$};
\node at (1+\dist,0) {$\C$};

\draw[arrows={-angle 90}, shorten >=9, shorten <=9]  (-2+\dist,0) -- (-2+0.5+\dist,1);
\draw[arrows={-angle 90}, shorten >=9, shorten <=9]  (-2+0.5+\dist,1) -- (-1+\dist,0);
\draw[arrows={-angle 90}, shorten >=9, shorten <=9]  (-1+\dist,0) -- (-1+0.5+\dist,1);
\draw[arrows={-angle 90}, shorten >=9, shorten <=9]  (-1+0.5+\dist,1) -- (0+\dist,0);
\draw[arrows={-angle 90}, shorten >=9, shorten <=9]  (0+\dist,0) -- (0+0.5+\dist,1);
\draw[arrows={-angle 90}, shorten >=9, shorten <=9]  (0+0.5+\dist,1) -- (1+\dist,0);

\node at (-2+0.5+\dist,1) {$\C$};
\node at (-1+0.5+\dist,1) {$\C$};
\node at (0+0.5+\dist,1) {$\C$};

\draw[arrows={-angle 90}, shorten >=9, shorten <=9]  (-2+0.5+\dist,1) -- (-2+1+\dist,2);
\draw[arrows={-angle 90}, shorten >=9, shorten <=9]  (-2+1+\dist,2) -- (-1+0.5+\dist,1);
\draw[arrows={-angle 90}, shorten >=9, shorten <=9]  (-1+0.5+\dist,1) -- (-1+1+\dist,2);
\draw[arrows={-angle 90}, shorten >=9, shorten <=9]  (-1+1+\dist,2) -- (0+0.5+\dist,1);

\node at (-1+\dist,2) {$\C$};
\node at (0+\dist,2) {$\C$};

\draw[arrows={-angle 90}, shorten >=9, shorten <=9]  (-1+\dist,2) -- (-1+0.5+\dist,3);
\draw[arrows={-angle 90}, shorten >=9, shorten <=9]  (-1+0.5+\dist,3) -- (-0+\dist,2);

\node at (-1+0.5+\dist,3) {$\C$};

}

\foreach \dist in {2}{
\draw[arrows={-angle 90}, shorten >=5, shorten <=5,->>]  (-0.5+\dist,2) -- (0.5+\dist,2);
}

\foreach \dist in {5}{
\node at (0-0.5+\dist,4) {$I_{i-5,4}$};
\node at (-2+\dist,0) {$\C$};
\draw[arrows={-angle 90}, shorten >=9, shorten <=9]  (-2+\dist,0) -- (-2+0.5+\dist,1);
\node at (-2+0.5+\dist,1) {$\C$};
\draw[arrows={-angle 90}, shorten >=9, shorten <=9]  (-2+0.5+\dist,1) -- (-2+1+\dist,2);
\node at (-1+\dist,2) {$\C$};
\draw[arrows={-angle 90}, shorten >=9, shorten <=9]  (-1+\dist,2) -- (-1+0.5+\dist,3);
\node at (-1+0.5+\dist,3) {$\C$};

}

\end{tikzpicture}
\]
\end{ex}
\subsection{The Restriction Functor}
For each $W \in \mr{rep}_\C(\hat{\Delta}_{n,N},\hat{\mr{I}}_{n,N})$ we define the representation $\mr{res} \, W \in \mr{rep}_\C(\Delta_n,\mr{I}_N)$ as
\[ \mr{res}  \, W := \Big( (W^{(i,1)})_{i \in \Z_n}, (W_{\beta_{i,2}} \circ W_{\alpha_{i,1}})_{i \in \Z_n} \Big).\] 
This induces maps $\mr{res}_{V,W} : \mr{Hom}_{\hat{\Delta}_{n,N}}(V,W) \to \mr{Hom}_{\Delta_n}(\mr{res}  \,V,\mr{res}  \,W)$, by forgetting the components of the morphisms at the vertices $(i,k)$ with $k\geq 2$. Hence we obtain a functor $\mr{res} :  \mr{rep}_\C(\hat{\Delta}_{n,N},\hat{\mr{I}}_{n,N}) \to  \mr{rep}_\C(\Delta_n,\mr{I}_N)$. 
The proof of the following proposition is immediate by the construction of $\Lambda$ and $\mr{res}$.
\begin{prop}
$\mr{res} \circ \Lambda (M) = M$ holds for all $M \in  \mr{rep}_\C(\Delta_n,\mr{I}_N)$.
\end{prop}
\subsection{The Desingularization Map}
In this subsection we provide the construction of the desingularization map, again closely following \cite{CFR2013}.  An example is given in Section~\ref{sec:example}.
\begin{defi}
An isomorphis class $[N]$ of $\Delta_n$-representations is called a \f{generic subrepresentation type} of $M \in \mr{rep}_\C(\Delta_n,\mr{I}_N)$ to dimension vector $\mb{e}$, if the stratum $\mathcal{S}_{[N]}$ is open in $\mr{Gr}_\mb{e}(M)$.
The set of generic subrepresentation types is denoted by $\mr{gsub}_\mb{e}(M)$. 
\end{defi}


\begin{rem}
By construction, the closure of the stratum $\mathcal{S}_{[N]}$ for some $[N] \in \mr{gsub}_\mb{e}(M)$ is an irreducible component of $\mr{Gr}_\mb{e}(M)$, and all irreducible components are obtained in this way.
\end{rem}

\begin{rem}
In general, there is no explicit description of the $\mr{gsub}_\mb{e}(M)$. But if the indecomposable summands of $M$ are all of length $\omega n$ for $n,\omega \in \N$, we can apply Lemma~\ref{lma:irr_comp}.
\end{rem}


\begin{ex}
Let $n=3$, $N=2$ and consider the quiver Grassmannian for $M=U_1(2)^2\oplus U_2(2)^3\oplus U_3(2)$ and $\mb{e}=(1,2,3)$. It has eight isomorphism classes of subrepresentations but only two irreducible components. Namely $\overline{\mathcal{S}_{[N_{1,2}]}}$ for $N_1=U_2(2)^2\oplus U_3(2)$ and $N_2=U_1(2)\oplus U_2(2)\oplus U_3(1)^2$. The stratum of $N_1$ is $7$-dimensional whereas the stratum of $N_2$ is only $5$-dimensional. 
\end{ex}

For $[N] \in \mr{gsub}_\mb{e}(M)$ we define a map from a quiver Grassmannian of the cylinder quiver to a quiver Grassmannian of the cycle quiver
\[ \pi_N : \mr{Gr}_{\bdim \hat{N} } (\hat{M}) \to \mr{Gr}_\mb{e}(M) \]
by $\pi_N(U) := \mr{res} \, U$ for all $U \in \mr{Gr}_{\bdim \hat{N} } (\hat{M})$. 

\begin{prop}\label{prop:injective-over-strata}
For each $[N] \in \mr{gsub}_\mb{e}(M)$ the map
\[ \pi_N : \mr{Gr}_{\bdim \hat{N} } (\hat{M}) \to \mr{Gr}_\mb{e}(M) \]
is injective over $\mathcal{S}_{[N]}$.
\end{prop}

\begin{proof}
Let $U \in \mathcal{S}_{[N]} \subseteq  \mr{Gr}_\mb{e}(M)$, then $\bdim \hat{U} = \bdim \hat{N}$ and
\[ \pi^{-1}_N(U) = \Big\{ V \in \mr{Gr}_{\bdim \hat{N} } (\hat{M}) \ : \ V^{(i,1)} = U^{(i)} \fa i \in \Z_n\Big\}. \]
In particular $\hat{U}$ is contained in $\pi^{-1}_N(U) \subset \mr{Gr}_{\bdim \hat{N} } (\hat{M})$. It remains to show that $\pi^{-1}_N(U) = \{ \hat{U} \}$: By construction of $\mr{res}$ and $\Lambda$ it follows that $\hat{U}^{(i,2)} \subset V^{(i,2)}$ and $\dim_\C \hat{U}^{(i,2)}  = \dim_\C V^{(i,2)}$ holds for all $V  \in \pi^{-1}_N(U)$ since $U$ and $N$ are isomorphic. This implies that $\hat{U}^{(i,2)} = V^{(i,2)}$ holds for all $i \in \Z_n$. Inductively, it follows that $V = \hat{U}$.
\end{proof}

\begin{prop}\label{prop:fibre-of-pi}
For each $[N] \in \mr{gsub}_\mb{e}(M)$ the fibre of $\pi_N$ over $U \in \mr{Gr}_\mb{e}(M) $ is 
\[ \pi_N^{-1}(U) = \mathcal{F}_U:= \big\{ F \in \mr{Gr}_{\bdim \hat{N} } (\hat{M}) \ \big\vert \ \hat{U} \subseteq F \big\} \cong  \mr{Gr}_{\bdim \hat{N} -\bdim \hat{U} } \big(\hat{M}/\hat{U} \big). \]
\end{prop}

\begin{proof}
Observe that $\bdim U = \bdim N$, so that $\dim_\C \hat{U}^{(i,1)} =\dim_\C \hat{N}^{(i,1)}$ for all $i\in \Z_n$ and the first non-trivial choice of a subspace $F^{(i,k)}$ is over vertices $(i,k)$ with $k \geq 2$. The inclusion $\mathcal{F}_U \subseteq \pi_N^{-1}(U)$ holds since $\pi_N(F)=U$ is clear by definition of $\mathcal{F}_U$ and the construction of the restriction functor. The other inclusion follows since every point $V$ of the fibre $\pi_N^{-1}(U)$ has to contain the vector spaces of $\hat{U}$ in its vector spaces $V^{(i,k)}$ over each vertex $(i,k)$ of $\hat{\Delta}_{n,N}$, in order to map to $U$. The isomorphism between $\mathcal{F}_U$ and the quiver Grassmannian is a direct consequence of the explicit description of the fibre. 
\end{proof}

We are now ready to state the main result of the paper, which is proved after the next proposition.
\begin{trm}\label{trm:desing}
Let  $M \in  \mr{rep}_\C(\Delta_n,\mr{I}_N)$.
The map
\[ \pi := \bigsqcup_{[N] \in \mr{gsub}_\mb{e}(M)} \pi_N  :  \bigsqcup_{[N] \in \mr{gsub}_\mb{e}(M)} \mr{Gr}_{\bdim \hat{N} } (\hat{M}) \to  \mr{Gr}_\mb{e}(M)\]
is a  desingularization of $\mr{Gr}_\mb{e}(M)$.
\end{trm}

\begin{rem}
Using Proposition~\ref{prop:fibre-of-pi}, we can compute the fibre dimensions for the desingularization to examine whether it is small, in the spirit of \cite[Section 2]{FF13}. This is the case for the quiver Grassmannian $\mr{Gr}_2(M)$ from \cite[Example~3.13]{LaPu2021} where $Q=\Delta_1$ and $M=U_1(2)\oplus S_1^2$. In general,  desingularizations of quiver Grassmannians for the cycle are not small. It already fails for the loop quiver (i.e. $\Delta_1$) and the quiver Grassmannian $\mr{Gr}_2(N)$ where $N = U_1(2)^2$.
\end{rem}

For the proof of Theorem~\ref{trm:desing} we recollect the main properties of the maps $\pi_N$:

\begin{prop}\label{prop:properties-of-the-maps}
Let $M \in \mr{rep}_\C(\Delta_n,\mr{I}_N)$ and $[N] \in \mr{gsub}_\mb{e}(M)$. Then:
\begin{itemize}
\item[(i)] The variety $\mr{Gr}_{\bdim \hat{N}}(\hat{M})$ is smooth with irreducible equidimensional connected components.
\item[(ii)] The map $\pi_N$ is one-to-one over $\mathcal{S}_{[N]}$. 
\item[(iii)] The image of $\pi_N$ is closed in $\mr{Gr}_\mb{e}(M)$ and contains $\overline{\mathcal{S}_{[N]}}$.
\item[(iv)] The map $\pi_N$ is projective.
\end{itemize}
\end{prop}

\begin{proof}
By Theorem~\ref{trm:global-dim-two} and Lemma~\ref{lma:proj-inj-dim-one} we can apply  \cite[Proposition~7.1]{CFR2013} to each quiver Grassmannian $\mr{Gr}_{\bdim \hat{N}}(\hat{M})$ and obtain the properties stated in (i). Proposition~\ref{prop:injective-over-strata} is exactly part (ii). The remaining parts are proven analogous to \cite[Theorem~7.5]{CFR2013} since the functor $\Lambda$ is fully faithful by Proposition~\ref{prop:Lambda-in-fully-faithful}.
\end{proof}

\begin{proof}[Proof of Theorem~\ref{trm:desing}]
By \cite[Proposition~37]{CEFR2018}, we obtain that \[\mr{Gr}_{\bdim \hat{N} } (\hat{M}) = \overline{\mathcal{S}_{[\hat{N}]}}\] since $\hat{M}$ is rigid by Lemma~\ref{lma:proj-inj-dim-one}. With the properties of $\hat{\Delta}_{n,N}$-representations from Theorem~\ref{trm:global-dim-two} and Lemma~\ref{lma:proj-inj-dim-one}, the maps $\pi_N$ as in Proposition~\ref{prop:properties-of-the-maps} and  $\Lambda$ as in Proposition~\ref{prop:Lambda-in-fully-faithful}, the rest of the proof is the same as for \cite[Corollary~7.7]{CFR2013}.
\end{proof}
\begin{rem}
In particular, \cite[Proposition~37]{CEFR2018} proves the conjecture from \cite[Remark~7.8]{CFR2013},
about the irreducibility of $\mr{Gr}_{\mb{e}} (\hat{M})$ in \cite[Corollary~7.7]{CFR2013} for arbitrary representations $M$ of a Dynkin quiver. 
\end{rem}
The following result generalizes \cite[Theorem~7.10]{FLP21}.

\begin{trm}\label{trm:tower-of-grassmann-bundles}
For each $[N] \in \mr{gsub}_\mb{e}(M)$ the quiver Grassmannian $\mr{Gr}_{\bdim \hat{N} } (\hat{M}) $ is isomorphic to a tower of fibrations 
\[  \mr{Gr}_{\bdim \hat{N} } (\hat{M}) = X_1 \to X_{2} \to \dots \to X_N = \prod_{i \in \Z_n}  \mr{Gr}_{\hat{n}^{(i,N)}}\Big(\C^{\hat{m}^{(i,N)}}\Big) \]
where $\hat{\mb{n}} := \bdim \hat{N}$ and $\hat{\mb{m}} := \bdim \hat{M}$ and each map $X_k \to X_k+1$ for $k \in [N-1]$ is a fibration with fibre isomorphic to a product of ordinary Grassmannians of subspaces.
\end{trm}

\begin{proof}
Every point $U$ of the quiver Grassmannian $\mr{Gr}_{\bdim \hat{N} } (\hat{M}) $ is parameterized by a collection of subspaces $U^{(i,k)} \subseteq M^{(i,k)}$ for $i \in \Z_n$ and $k \in [N]$. In particular it is a point in 
\[ \mr{Gr}_{\hat{\mb{n}}}  \big(\C^{\hat{\mb{m}}}\big) := \prod_{i \in \Z_n} \prod_{k \in [N]} \mr{Gr}_{\hat{n}^{(i,k)}}\Big(\C^{\hat{m}^{(i,k)}}\Big). \]
Define $X_k$ as the image of $\mr{Gr}_{\bdim \hat{N} } (\hat{M})$ in the variety  $\mr{Gr}_{\hat{\mb{n}}}({\hat{\mb{m}}}) ^{(k)}$ which is defined analogous to  $\mr{Gr}_{\hat{\mb{n}}}\big(\C^{\hat{\mb{m}}}\big) $, with the only difference that the second product  runs over $\{k,k+1,\dots,N\}$ instead of $[N]$. Hence $\mr{Gr}_{\bdim \hat{N} } (\hat{M}) = X_1 $ follows  by construction.

We proceed by decreasing induction starting from $k=N$. Every point in the the product of Grassmannians of subspaces $\mr{Gr}_{\hat{\mb{n}}}({\hat{\mb{m}}})^{(N)}$ can be extended to an element of $\mr{Gr}_{\bdim \hat{N} } (\hat{M})$ since the upper vector spaces of an element in the quiver Grassmannian are not related. This implies $X_N = \mr{Gr}_{\hat{\mb{n}}}({\hat{\mb{m}}})^{(N)}$ as desired. 

Now assume that the vector spaces $U^{(i,k')}$ are fixed for all $i \in \Z_n$ and $k' > k$. Since $U$ has to be contained in the quiver Grassmannian it has to satisfy the relations
\[ \hat{M}_{\alpha_{i+1,k}} \circ \hat{M}_{\beta_{i,k+1}} U^{(i,k+1)} \subseteq U^{(i+1,k+1)} \quad \fa i \in\Z_n .\] 
Hence the next layer of vector spaces $U^{(i,k)}$ requires 
\[ \hat{M}_{\beta_{i-1,k+1}} U^{(i-1,k+1)} \subseteq U^{(i,k)} \quad \mr{and} \quad \hat{M}_{\alpha_{i,k}}  U^{(i,k)}  \subseteq  U^{(i,k+1)} \quad \fa i \in\Z_n . \] 
This is equivalent to the choice of a point in the Grassmannian 
\[ \mr{Gr}_{\hat{n}^{(i,k)}-\hat{n}^{(i-1,k+1)}}\Big(U^{(i,k)}/\hat{M}_{\beta_{i-1,k+1}} U^{(i-1,k+1)}\Big) \]
because every map $\hat{M}_{\alpha_{i,k}}$ is a projection where the last $\hat{m}^{(i,k)}-\hat{m}^{(i,k+1)}$ coordinates are sent to zero and each $\hat{M}_{\beta_{i,k}}$ is an inclusion.
\end{proof}

\begin{rem}
The explicit description of the desingularization in Theorem~\ref{trm:tower-of-grassmann-bundles} allows to construct a cellular decomposition of $\mr{Gr}_{\bdim \hat{N} } (\hat{M})$ (c.f. Theorem~\ref{trm:cellular-decomp-desing}). In particular, it implies that $\mr{Gr}_{\bdim \hat{N} } (\hat{M}) $ is smooth. 
\end{rem}
\section{Torus Equivariant Cohomology and Equivariant Euler Classes}\label{sec:equivariant-cohomology}
In this section we briefly recall definitions and constructions concerning torus actions on quiver Grassmannians, torus equivariant cohomology and torus equivariant Euler classes. More details on the general theory is found in \cite{Arabia98,Brion1997,GKM1998,Gonzales2014}. The application to quiver Grassmannians is introduced in \cite{LaPu2020,LaPu2021}.
In Section~\ref{sec:application} we provide examples and apply our desingularizations to the computation of equivariant cohomology of quiver Grassmannians for the equioriented cycle.
\subsection{Moment Graph and Torus equivariant Cohomology }
Let $X$ be a projective algebraic variety over $\C$. The action of an algebraic torus $T \cong (\C^*)^r$ on $X$ is \f{skeletal} if the number of $T$-fixed points and the number of one-dimensional $T$-orbits in $X$ is finite. We call a cocharacter $\chi \in \mathfrak{X}_*(T)$ \f{generic} for the $T$-action on $X$ if $X^T = X^{\chi(\C^*)}$. By $\mathfrak{X}^*(T)$ we denote the character lattice of $T$. The $T$-equivariant cohomology of $X$ with rational coefficients is denoted by $H_T^\bullet(X)$. 
\begin{defi} The pair $(X,T)$ is a \f{GKM-variety} if the $T$-action on $X$ is skeletal and the rational cohomology of $X$ vanishes in odd degrees.
\end{defi}
\begin{rem}
By \cite[Lemma~2]{Brion2000} this is equivalent to \cite[Definition~1.4]{LaPu2020}. 
\end{rem}
The closure $\overline{E}$ of every one-dimensional $T$-orbit $E$ in a projective GKM-variety admits an $T$-equivariant isomorphism to $\C\mathbb{P}^1$. Thus each one-dimensional $T$-orbit connects two distinct $T$-fixed points of $X$.
\begin{defi}\label{def:moment-graph}Let $(X,T)$ be a GKM-variety, and let $\chi\in \mathfrak{X}_*(T)$ be a generic cocharacter. The corresponding \f{moment graph} $\mathcal{G}=\mathcal{G}(X,T, \chi)$ of a GKM-variety is given by the following data:
\begin{itemize}
\item[(MG0)] the $T$-fixed points as vertices, i.e.: $\mathcal{G}_0 = X^T$,
\item[(MG1)] the closures of one-dimensional $T$-orbits $\overline{E} = E \cup \{x,y\}$ as edges in $\mathcal{G}_1$, oriented from $x$ to $y$ if $\lim_{\lambda \to 0}\chi(\lambda).p=x$ for $p \in E$,
\item[(MG2)] every $\overline{E}$ is labelled by $\alpha_E \in \mathfrak{X}^*(T)$ describing the $T$-action on $E$. 
\end{itemize}
\end{defi}
\begin{trm}\label{thm:GKM}(\cite[Theorem 1.2.2]{GKM1998}) Let $(X, T)$ be a GKM-variety with moment graph $\mathcal{G}=\mathcal{G}(X,T,\chi)$ and set $R:=H_T^\bullet(\textrm{pt})$. Then 
\[
H_T^\bullet(X) \cong \left\{(f_x)\in\bigoplus_{x\in \mathcal{G}_0}R \ \Big| \
\begin{array}{c}
f_{x_E}-f_{y_E}\in \alpha_{E} R\\
 \hbox{ for any }\overline{E}=E\cup\{x_E, y_E\}\in\mathcal{G}_1\end{array}
\right\}.\]  
\end{trm}
\begin{rem}\label{rem:torus-characters-part-i}
The characters from (MG2) are only unique up to a sign. This sign does not play a role in Theorem \ref{thm:GKM}. Hence we can fix our favourite convention.
\end{rem}
\subsection{BB-filterable varieties}
In this subsection we describe a class of varieties which admit an explicit formula for the computation of their equivariant cohomology. Let $X$ be a $\C^*$-variety. By $X^{\C^*}$ we denote its fixed point set and $X_1, \ldots, X_m$ denote the connected components of $X^{\C^*}$.
This induces a decomposition
\begin{equation}\label{eqn:BBdecomposition}
X=\bigcup_{i\in [m]} W_i, \quad\hbox{ with }\quad   W_i := \left\{ x \in X \mid \lim_{z \to 0} z.x \in X_i \right\},
\end{equation}
where $W_i$ is called attracting set of $X_i$. Since decompositions of this type were first studied by Bialynicki-Birula in \cite{Birula1973}, we call it a \f{BB-decomposition}. 
\begin{defi}We say that $W_i$ from \eqref{eqn:BBdecomposition} is a \f{rational cell} if it is rationally smooth at all $w\in W_i$. This in turn holds if 
\[H^{2\textrm{dim}_{\C}(W_i)}(W_i, W_i\setminus \{w\})\simeq \Q  \quad\hbox{ and  }\quad  H^m(W_i,W_i\setminus\{w\})=0 \]
for any $m\neq 2 \textrm{dim}_{\C}(W_i)$ (cf. \cite[p.292, Definition~3.4]{Gonzales2014}). 
\end{defi}
\begin{defi}
\label{def:BB-filterable}
A projective $T$-variety $X$ is \f{BB-filterable} if: 
\begin{enumerate}
\item[(BB1)] the fixed point set $X^T$ is finite,
\item[(BB2)] there exists a generic cocharacter $\chi: \C^* \rightarrow T$, i.e. $X^{\chi(\C^*)} = X^T$, such that 
 the associated BB-decomposition consists of rational cells. 
\end{enumerate}
\end{defi}
 \begin{trm}
\label{trm:t-stable_filtration-general-setting}(cf. \cite[Theorem~1.15]{LaPu2020}) Let $X$ be a BB-filterable projective $T$-variety. Then: 
\begin{enumerate}
\item $X$ admits a filtration into $T$-stable closed subvarieties $Z_i$ such that
\[ \emptyset = Z_0 \subset Z_1 \subset \dots \subset Z_{m-1} \subset Z_m = X. \]
\item Each $W_i = Z_i \setminus Z_{i-1}$ is a rational cell, for all $i \in [m]$.
\item The singular rational cohomology of $Z_i$ vanishes in odd degrees, for $i \in [m]$. 
\item If, additionally, the $T$-action on $X$ is skeletal, each $Z_i$ is a GKM-variety.
\end{enumerate}
\end{trm}
\subsection{Euler Classes and Cohomology Module Bases}
For the precise definition of Euler Classes we refer the reader to \cite[Section~2.2.1]{Arabia98}. Instead we give three properties which are enough to determine the equivariant Euler classes in our setting. 
\begin{lem}\label{lma:euler-class-along-resolution}(cf. \cite[Corollary~15, Lemma~16, Theorem~18]{Brion1997}) Let $Y$ be a $T$-variety and $y\in Y^T$. 
 \begin{enumerate}
 \item If $Y$ is smooth at $y$ then $\mr{Eu}_T(y,Y)=(-1)^{\mr{dim}(Y)} \mr{det} \, T_yY$, where $\mr{det} \, T_yY$ is the product of the characters by which $T$ acts on the tangent space $T_yY$.
 \item If $Y$ is rationally smooth at $y$ then $\mr{Eu}_T(y,Y)=z\cdot \mr{det} \, T_yY$, for some $z\in\Q\setminus\{0\}$. 
\item If $\pi: Y\to X$ is a $T$-equivariant resolution of singularities and $|Y^T|<\infty$, then
\[ \mr{Eu}_T(x,X)^{-1} = \sum_{ y \in Y^T, \pi(y)=x} \mr{Eu}_T(y,Y)^{-1}.\]
 \end{enumerate}
 \end{lem}
\begin{rem} 
Lemma~\ref{lma:euler-class-along-resolution} differs from \cite{Brion1997} by using Euler classes instead of equivariant multiplicities which are inverse to each other up to a sign.
 \end{rem}
 \begin{defi}(cf. \cite[Lemma~6.7]{Gonzales2014})\label{def:localIndex} Let $X^T=\{x_1, \ldots, x_m\}$. For $i \in [m]$, the local index of $f\in H_T^\bullet(X)$ at $x_i \in X^T$ is 
\begin{equation*}
I_i(f) = \sum_{\substack{j \in [m] \ : \ \\x_j\in Z_i}} \frac{f_{x_j}}{\mr{Eu}_T(x_j,Z_i) \ \ }.
\end{equation*}
\end{defi}
The next theorem gives an explicit formula to compute a basis for $H_T^\bullet(X)$ as free module over $H_T^\bullet(pt)$. Observe that everything depends on the order of the fixed points which is in general not unique.
\begin{trm}
\label{trm:cohomology-generators-general-setting}(cf. \cite[Theorem~2.12]{LaPu2020}) Let $(X,T)$ be a BB-filterable GKM-variety with filtration 
\[ \emptyset=Z_0\subset Z_1\subset\ldots\subset Z_m=X\] 
as in Theorem~\ref{trm:t-stable_filtration-general-setting}. Let $X^T=\{x_1, \ldots, x_m\}$ with $x_i\in W_i=Z_i\setminus Z_{i-1}$. There exists a unique basis $\{\theta^{(i)}\}_{i \in [m]}$ of $H_T^\bullet(X)$ as free module over $H_T^\bullet(pt)$, such that for any $i \in [m]$ the following properties hold:
\begin{enumerate}
\item  $\theta^{(i)}_{x_j}=0$ for all $j < i$,
\item  $\theta^{(i)}_{x_i}= \mr{Eu}_T(x_i,Z_i)$,
\item  $I_j(\theta^{(i)}) = 0$  for all $j \neq i$.
 \end{enumerate}
\end{trm}
\begin{rem}
Observe that (1) and (2) imply $I_i(\theta^{(i)})= 1$ by Definition~\ref{def:localIndex}.
\end{rem}
\subsection{Torus Action on Cyclic Quiver Grassmannians}
We briefly recall torus actions on quiver Grassmannians for the equioriented cycle (c.f. \cite[Section 5]{LaPu2020}).
\begin{rem}
From now on we assume the choice of a basis $B$ of $M$ such that the connected components of $Q(M,B)$ are in bijection with the indecomposable direct summands of $M$. Such a choice is always possible by \cite[Theorem~1.11]{Kirillov2016}.
\end{rem}
A \f{grading} of $M \in \mr{rep}_\C(\Delta_n)$ with respect to a fixed basis is a map $\mr{wt} : B \to \Z^B$. This induces an action of $\lambda \in \C^*$ by 
\[ \lambda.b := \lambda^{\mr{wt}(b)}\cdot b. \]
\begin{rem}
Combining several weight functions $\mr{wt}_1, \dots, \mr{wt}_D : B \to \Z^B$, we can define the action of $\lambda = ( \lambda_j)_{j \in [D]} \in (\C^*)^D$ as 
\[ \lambda.b := \prod_{j \in [D]} \lambda_j^{\mr{wt}_j(b)} \cdot b =  \lambda_1^{\mr{wt}_1(b)} \cdot \ldots  \cdot \lambda_D^{\mr{wt}_D(b)} \cdot b. \] 
Observe that this action extends to the quiver Grassmannian $\mr{Gr}_\mb{e}(M)$ only under some additional assumptions about the grading (c.f. \cite[Lemma~5.12]{LaPu2020}). 
\end{rem}
\begin{trm}
\label{trm:t-stable_filtration-equi-cycle} (c.f. \cite[Theorem~6.6]{LaPu2020})
Let $M$ be a nilpotent representation of $\Delta_n$ with $d$-many indecomposable direct summands, and let $\mb{e} \leq \bdim M$ be such that $\mr{Gr}_\mb{e}(M)$ is non-empty. Let $T := (\C^{*})^{d+1}$ act on $\mr{Gr}_\mb{e}(M)$ as in \cite[Lemma~5.12]{LaPu2020}. Then $(\mr{Gr}_\mb{e}(M),T)$ is a projective BB-filterable GKM-variety.
\end{trm}
\begin{rem}
If the desingularizations constructed in Theorem\ref{trm:desing} are $T$-equivariant, this theorem implies that we can compute the $T$-equivariant cohomology of all quiver Grassmannians for nilpotent representations of $\Delta_n$, using Theorem~\ref{trm:cohomology-generators-general-setting}.
\end{rem}
From now on we assume that $T := (\C^*)^{d+1}$ acts on $\mr{Gr}_\mb{e}(M)$ as in \cite[Lemma~5.12]{LaPu2020}. Here $d$ is the number of connected components in $Q(M,B)$ and the additional parameter comes from cyclic symmetry. The weight functions of the action are defined implicitly by the formula used in \cite[Section~5.2]{LaPu2020}.
\section{Torus Equivariant Desingularization and Application}\label{sec:application}
 In this section we apply the methods from the previous section to compute Euler classes at singular points and torus equivariant cohomology of quiver Grassmannians for the equioriented cycle using their desingularizations as constructed in Theorem~\ref{trm:desing}. It remains to show that these desingularizations are torus equivrariant.
\subsection{Torus Action on the Desingularization}
Let $M \in \mr{rep}_\C(\Delta_n,\mr{I}_N)$ be nilpotent with $d$-many indecomposable direct summands, and let $T := (\C^{*})^{d+1}$ act on $\mr{Gr}_\mb{e}(M)$ as in \cite[Lemma~5.12]{LaPu2020}.
\begin{rem}\label{rem:T-action-desing}
A choice of basis $B$ of $M \in  \mr{rep}_\C(\Delta_n,\mr{I}_N)$ induces a basis $\hat{B}$ of $\hat{M} \in  \mr{rep}_\C(\hat{\Delta}_{n,N},\hat{\mr{I}}_{n,N})$ such that the connected components of $Q(\hat{M},\hat{B})$ are in bijection with the images of the indecomposable summands of $M$. In particular the basis $\hat{B}^{(i,k)}$ over the vertex $(i,k)$ of $\hat{\Delta}_{n,N}$ is a subset in the basis $B^{(i+k-1)}$ of cardinality $m_{i+k-1} -c$ where $c$ is the corank of the map $M_{\alpha_{i+k-2}} \circ \dots \circ M_{\alpha_{i}}$ if $k \geq 2$ and $\hat{B}^{(i,k)} = B^{(i)}$ for $k=1$. This allows us to extend the $T$-action to the vector spaces of $\hat{M}$ by extending the weight functions according to the inclusions of the basis described above. In other words, all basis vectors of $\hat{B}$ which have the same image in $B$ get the same weight. 
\end{rem}
\begin{prop}\label{prop:T-action-extends-to-desing}
The $T$-action on the vector spaces of $\hat{M}$ as defined in Remark~\ref{rem:T-action-desing} extends to every quiver Grassmannian $\mr{Gr}_\mb{k}(\hat{M})$.
\end{prop}
\begin{proof}
We have to show that the $T$-action is compatible with the maps of the quiver representation $\hat{M}$. By construction of the action and the representation $\hat{M}$, this follows immediately from the compatibility of the $T$-action (on the vector spaces of $M$) with the maps of $M$ as shown in \cite[Lemma~5.12]{LaPu2020}.
\end{proof}
\begin{lma}\label{lem:T-equiv-desing}
The  desingularization of Theorem~\ref{trm:desing} is $T$-equivariant.
\end{lma}
\begin{proof}
With Proposition~\ref{prop:T-action-extends-to-desing}, the statement follows immediately from the construction of the grading as in Remark~\ref{rem:T-action-desing} together with the description of the desingularization in Theorem~\ref{trm:desing}.
\end{proof}
\begin{rem}
The $T$-equivariance of the  desingularization allows us to use \cite[Lemma~2.1.(3)]{LaPu2020} for the computation of equivariant Euler classes at the singular points of $\mr{Gr}_\mb{e}(M)$. This allows us to apply \cite[Theorem~2.12]{LaPu2020} about the construction of a basis for the $T$-equivariant cohomology to all quiver Grassmannians for nilpotent representations of the cycle.
\end{rem}
\subsection{Cellular Decomposition of the Desingularization}
\begin{trm}\label{trm:cellular-decomp-desing}
For $[N] \in \mr{gsub}_{\mb{e}}(M)$ the $T$-fixed points of $\mr{Gr}_{\bdim \hat{N} } (\hat{M})$ are exactly the preimages of the $T$-fixed points of $\overline{\mathcal{S}_{[N]}} \subset \mr{Gr}_\mb{e}(M)$ under $\pi_N$. The $\C^*$-attracting sets of these points provide a cellular decomposition of $\mr{Gr}_{\bdim \hat{N} } (\hat{M})$.
\end{trm}
\begin{proof}
The $T$-equivariance of $\pi_N$ from Lemma~\ref{lem:T-equiv-desing} gives the desired description of the fixed points. Now we prove that the $\C^*$-attracting sets of these fixed points from the BB-decomposition are cells. By \cite[Lemma~4.12]{Carrell02}, they provide an $\alpha$-partition, i.e. there exists a total order of the fixed points 
$ \mr{Gr}_{\bdim \hat{N} } (\hat{M})^{\mathbb{C}^*} = \{p_1, \dots, p_r \}$ such that \( \bigsqcup_{j=1}^s W_i \)
is closed in $\mr{Gr}_{\bdim \hat{N} } (\hat{M})$ for all $s \in [r]$. It remains to show that they are isomorphic to affine spaces. This is induced by the cellular decomposition of $\mr{Gr}_\mb{e}(M)$ and the $T$-equivariance of the desingularization:

Let $p \in \mr{Gr}_{\bdim \hat{N} } (\hat{M})$ be a $T$-fixed point. The vector space $p^{(i,k)}$ over the vertex $(i,k)$ of $\hat{\Delta}_{n,N}$ is a point in the Grassmannian of subspaces $\mr{Gr}_{\hat{n}^{(i,k)}}\big(\C^{\hat{m}^{(i,k)}}\big)$. By construction of the $\C^*$-action (as in Remark~\ref{rem:T-action-desing}), the attracting set of $p^{(i,k)}$ in $\mr{Gr}_{\hat{n}^{(i,k)}}\big(\C^{\hat{m}^{(i,k)}}\big)$ is a cell. The attracting set of $p$ in the whole quiver Grassmannian is the intersection of these cells along the maps of $\hat{M}$. We proceed by induction on $k$. For $k=1$ there is nothing to show because there are no maps between the vector spaces. If $k=2$, we have the original vector spaces of the representation $M$ and one additional layer of subspaces therein. The relations between the coordinates in the attracting sets are the same as for $\mr{Gr}_\mb{e}(M)$. Hence they are cells by \cite[Theorem~5.7]{LaPu2020}. The maps of $\hat{M}$ along the arrows $\beta_{i,k}$ of $\hat{\Delta}_{n,N}$ are inclusions and the maps along $\alpha_{i,k}$ are projections where the last $m^{(i,k)}-m^{(i,k+1)}$ coordinates are sent to zero (c.f. \cite[Proposition~4.8]{LaPu2020}). Thus we obtain that the intersecting relations for each $k\in [N]$ are of the form as described in \cite[Theorem~5.7]{LaPu2020}. This implies the desired isomorphisms to affine spaces. 
\end{proof}
\begin{rem}
In the setting that 
\[ M = \bigoplus_{i \in \Z_n} U_i(\omega n) \quad \mr{and} \quad \mb{e} = (\omega k,\dots,\omega k) \in \Z^n \]
it is possible to strengthen the results concerning the desingularization (c.f. \cite[\S~2.5 \& 2.6]{FLP2023}). Namely, $\mr{Gr}_\mb{e}(M)$ has $\binom{n}{k}$ explicitly described irreducible components (see Example~\ref{ex:generic-subrep-special-case}) and the cells of $\mr{Gr}_{\bdim \hat{N} } (\hat{M})$ are the strata of the corresponding $T$-fixed points.
\end{rem}
\subsection{Example}\label{sec:example}
Now, we provide an explicit example for the constructions from the previous sections. Let $M := U_1(4)\oplus U_2(2) \oplus U_2(2)$ be a $\Delta_2$-representation and fix the dimension vector $\mb{e}=(2,2)$. The quiver Grassmannian $\mr{Gr}_\mb{e}(M)$ has five strata (i.e. isomorphism classes of subrepresentations) with the representatives:
\[ V_1:=U_1(4), \quad V_2:=S_2\oplus U_1(3), \quad V_3:=U_2(2)\oplus U_2(2),\]
\[V_4:=U_1(2) \oplus U_2(2), \quad V_5:=S_1 \oplus S_2 \oplus U_2(2). \]
The stratum of $V_2$ is three-dimensional, the strata of $V_1$, $V_3$ and $V_4$ are two-dimensional and the stratum of $V_5$ is one dimensional. This is computed using \cite[Proposition~4.4]{Pue2020} and \cite[Lemma~2.4]{CFR2012}. 

Let the basis $B$ of $M$ be the union of the standard basis for each indecomposable summand of $M$. Then its coefficient quiver is 
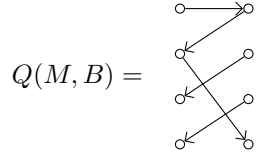
\begin{center}
\begin{tikzpicture}[scale=0.3]
\node at (-6,-2) {$Q(M,B) = $};
\draw[fill=white] (1.5,1) circle (.2);
\draw[fill=white] (-1.5,1) circle (.2);
\draw[fill=white] (1.5,-1) circle (.2);
\draw[fill=white] (-1.5,-1) circle (.2);
\draw[fill=white] (1.5,-3) circle (.2);
\draw[fill=white] (-1.5,-3) circle (.2);
\draw[fill=white] (1.5,-5) circle (.2);
\draw[fill=white] (-1.5,-5) circle (.2);

\draw[arrows={-angle 90}, shorten >=2, shorten <=2]  (-1.5,1) -- (1.5,1);
\draw[arrows={-angle 90}, shorten >=2, shorten <=2]  (1.5,1) -- (-1.5,-1);
\draw[arrows={-angle 90}, shorten >=2, shorten <=2]  (-1.5,-1) -- (1.5,-5);

\draw[arrows={-angle 90}, shorten >=2, shorten <=2]  (1.5,-1) -- (-1.5,-3);
\draw[arrows={-angle 90}, shorten >=2, shorten <=2]  (1.5,-3) -- (-1.5,-5);

\end{tikzpicture}
\end{center}
where the arrows from left to right have $\alpha: 1 \to 2$ as underlying arrow in $\Delta_2$. The arrows from right to left have $\beta: 2 \to 1$ as underlying arrow. We define the action of $\gamma:= (\gamma_0,\gamma_1,\gamma_2,\gamma_3) \in T := (\C^*)^{3+1}$ and $\lambda \in \C^*$ on $B$ as follows:
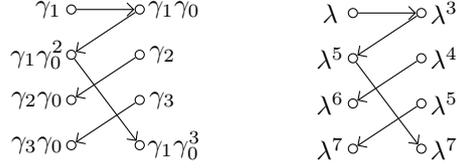
\begin{center}
\begin{tikzpicture}[scale=0.3]

\draw[fill=white] (1.5,1) circle (.2);\node at (3,1) {$\gamma_1\gamma_0$};
\draw[fill=white] (-1.5,1) circle (.2);\node at (-2.5,1) {$\gamma_1$};
\draw[fill=white] (1.5,-1) circle (.2);\node at (2.5,-1) {$\gamma_2$};
\draw[fill=white] (-1.5,-1) circle (.2);\node at (-3,-1) {$\gamma_1\gamma_0^2$};
\draw[fill=white] (1.5,-3) circle (.2);\node at (2.5,-3) {$\gamma_3$};
\draw[fill=white] (-1.5,-3) circle (.2);\node at (-3,-3) {$\gamma_2\gamma_0$};
\draw[fill=white] (1.5,-5) circle (.2);\node at (3,-5) {$\gamma_1\gamma_0^3$};
\draw[fill=white] (-1.5,-5) circle (.2);\node at (-3,-5) {$\gamma_3\gamma_0$};

\draw[arrows={-angle 90}, shorten >=2, shorten <=2]  (-1.5,1) -- (1.5,1);
\draw[arrows={-angle 90}, shorten >=2, shorten <=2]  (1.5,1) -- (-1.5,-1);
\draw[arrows={-angle 90}, shorten >=2, shorten <=2]  (-1.5,-1) -- (1.5,-5);

\draw[arrows={-angle 90}, shorten >=2, shorten <=2]  (1.5,-1) -- (-1.5,-3);
\draw[arrows={-angle 90}, shorten >=2, shorten <=2]  (1.5,-3) -- (-1.5,-5);

\end{tikzpicture}
$\quad \quad \quad$
\begin{tikzpicture}[scale=0.3]
\draw[fill=white] (1.5,1) circle (.2);\node at (2.5,1) {$\lambda^3$};
\draw[fill=white] (-1.5,1) circle (.2);\node at (-2.5,1) {$\lambda$};
\draw[fill=white] (1.5,-1) circle (.2);\node at (2.5,-1) {$\lambda^4$};
\draw[fill=white] (-1.5,-1) circle (.2);\node at (-2.5,-1) {$\lambda^5$};
\draw[fill=white] (1.5,-3) circle (.2);\node at (2.5,-3) {$\lambda^5$};
\draw[fill=white] (-1.5,-3) circle (.2);\node at (-2.5,-3) {$\lambda^6$};
\draw[fill=white] (1.5,-5) circle (.2);\node at (2.5,-5) {$\lambda^7$};
\draw[fill=white] (-1.5,-5) circle (.2);\node at (-2.5,-5) {$\lambda^7$};

\draw[arrows={-angle 90}, shorten >=2, shorten <=2]  (-1.5,1) -- (1.5,1);
\draw[arrows={-angle 90}, shorten >=2, shorten <=2]  (1.5,1) -- (-1.5,-1);
\draw[arrows={-angle 90}, shorten >=2, shorten <=2]  (-1.5,-1) -- (1.5,-5);

\draw[arrows={-angle 90}, shorten >=2, shorten <=2]  (1.5,-1) -- (-1.5,-3);
\draw[arrows={-angle 90}, shorten >=2, shorten <=2]  (1.5,-3) -- (-1.5,-5);

\end{tikzpicture}
\end{center}
These actions extend linearly to the vector spaces of $M$ and to the whole quiver Grassmannian by \cite[Lemma~5.12]{LaPu2020}. Moreover 
\[ \chi : \C^* \to T; \quad  \lambda \mapsto (\lambda^2,\lambda,\lambda,\lambda)\]
is a generic cocharacter by \cite[Theorem~5.14]{LaPu2020}. We apply \cite[Theorem~1]{Cerulli2011} to compute the fixed points of both actions:
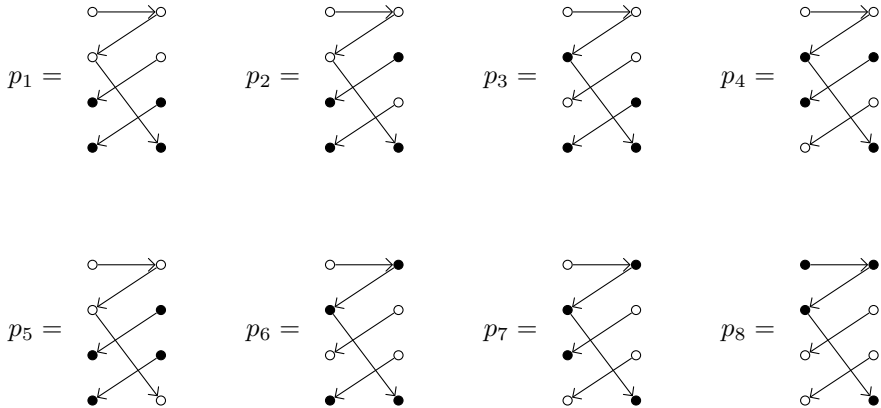
\begin{center}
\begin{tikzpicture}[scale=0.3]
\node at (-4,-2) {$p_1 = $};
\draw[fill=white] (1.5,1) circle (.2);
\draw[fill=white] (-1.5,1) circle (.2);
\draw[fill=white] (1.5,-1) circle (.2);
\draw[fill=white] (-1.5,-1) circle (.2);
\draw[fill=black] (1.5,-3) circle (.2);
\draw[fill=black] (-1.5,-3) circle (.2);
\draw[fill=black] (1.5,-5) circle (.2);
\draw[fill=black] (-1.5,-5) circle (.2);

\draw[arrows={-angle 90}, shorten >=2, shorten <=2]  (-1.5,1) -- (1.5,1);
\draw[arrows={-angle 90}, shorten >=2, shorten <=2]  (1.5,1) -- (-1.5,-1);
\draw[arrows={-angle 90}, shorten >=2, shorten <=2]  (-1.5,-1) -- (1.5,-5);

\draw[arrows={-angle 90}, shorten >=2, shorten <=2]  (1.5,-1) -- (-1.5,-3);
\draw[arrows={-angle 90}, shorten >=2, shorten <=2]  (1.5,-3) -- (-1.5,-5);

\end{tikzpicture}
$\quad \quad$
\begin{tikzpicture}[scale=0.3]
\node at (-4,-2) {$p_2 = $};
\draw[fill=white] (1.5,1) circle (.2);
\draw[fill=white] (-1.5,1) circle (.2);
\draw[fill=black] (1.5,-1) circle (.2);
\draw[fill=white] (-1.5,-1) circle (.2);
\draw[fill=white] (1.5,-3) circle (.2);
\draw[fill=black] (-1.5,-3) circle (.2);
\draw[fill=black] (1.5,-5) circle (.2);
\draw[fill=black] (-1.5,-5) circle (.2);

\draw[arrows={-angle 90}, shorten >=2, shorten <=2]  (-1.5,1) -- (1.5,1);
\draw[arrows={-angle 90}, shorten >=2, shorten <=2]  (1.5,1) -- (-1.5,-1);
\draw[arrows={-angle 90}, shorten >=2, shorten <=2]  (-1.5,-1) -- (1.5,-5);

\draw[arrows={-angle 90}, shorten >=2, shorten <=2]  (1.5,-1) -- (-1.5,-3);
\draw[arrows={-angle 90}, shorten >=2, shorten <=2]  (1.5,-3) -- (-1.5,-5);

\end{tikzpicture}
$\quad \quad$
\begin{tikzpicture}[scale=0.3]
\node at (-4,-2) {$p_3 = $};
\draw[fill=white] (1.5,1) circle (.2);
\draw[fill=white] (-1.5,1) circle (.2);
\draw[fill=white] (1.5,-1) circle (.2);
\draw[fill=black] (-1.5,-1) circle (.2);
\draw[fill=black] (1.5,-3) circle (.2);
\draw[fill=white] (-1.5,-3) circle (.2);
\draw[fill=black] (1.5,-5) circle (.2);
\draw[fill=black] (-1.5,-5) circle (.2);

\draw[arrows={-angle 90}, shorten >=2, shorten <=2]  (-1.5,1) -- (1.5,1);
\draw[arrows={-angle 90}, shorten >=2, shorten <=2]  (1.5,1) -- (-1.5,-1);
\draw[arrows={-angle 90}, shorten >=2, shorten <=2]  (-1.5,-1) -- (1.5,-5);

\draw[arrows={-angle 90}, shorten >=2, shorten <=2]  (1.5,-1) -- (-1.5,-3);
\draw[arrows={-angle 90}, shorten >=2, shorten <=2]  (1.5,-3) -- (-1.5,-5);

\end{tikzpicture}
$\quad \quad$
\begin{tikzpicture}[scale=0.3]
\node at (-4,-2) {$p_4 = $};
\draw[fill=white] (1.5,1) circle (.2);
\draw[fill=white] (-1.5,1) circle (.2);
\draw[fill=black] (1.5,-1) circle (.2);
\draw[fill=black] (-1.5,-1) circle (.2);
\draw[fill=white] (1.5,-3) circle (.2);
\draw[fill=black] (-1.5,-3) circle (.2);
\draw[fill=black] (1.5,-5) circle (.2);
\draw[fill=white] (-1.5,-5) circle (.2);

\draw[arrows={-angle 90}, shorten >=2, shorten <=2]  (-1.5,1) -- (1.5,1);
\draw[arrows={-angle 90}, shorten >=2, shorten <=2]  (1.5,1) -- (-1.5,-1);
\draw[arrows={-angle 90}, shorten >=2, shorten <=2]  (-1.5,-1) -- (1.5,-5);

\draw[arrows={-angle 90}, shorten >=2, shorten <=2]  (1.5,-1) -- (-1.5,-3);
\draw[arrows={-angle 90}, shorten >=2, shorten <=2]  (1.5,-3) -- (-1.5,-5);

\end{tikzpicture}
\end{center}
$ $
\begin{center}

\begin{tikzpicture}[scale=0.3]
\node at (-4,-2) {$p_5 = $};
\draw[fill=white] (1.5,1) circle (.2);
\draw[fill=white] (-1.5,1) circle (.2);
\draw[fill=black] (1.5,-1) circle (.2);
\draw[fill=white] (-1.5,-1) circle (.2);
\draw[fill=black] (1.5,-3) circle (.2);
\draw[fill=black] (-1.5,-3) circle (.2);
\draw[fill=white] (1.5,-5) circle (.2);
\draw[fill=black] (-1.5,-5) circle (.2);

\draw[arrows={-angle 90}, shorten >=2, shorten <=2]  (-1.5,1) -- (1.5,1);
\draw[arrows={-angle 90}, shorten >=2, shorten <=2]  (1.5,1) -- (-1.5,-1);
\draw[arrows={-angle 90}, shorten >=2, shorten <=2]  (-1.5,-1) -- (1.5,-5);

\draw[arrows={-angle 90}, shorten >=2, shorten <=2]  (1.5,-1) -- (-1.5,-3);
\draw[arrows={-angle 90}, shorten >=2, shorten <=2]  (1.5,-3) -- (-1.5,-5);

\end{tikzpicture}
$\quad \quad$
\begin{tikzpicture}[scale=0.3]
\node at (-4,-2) {$p_6 = $};
\draw[fill=black] (1.5,1) circle (.2);
\draw[fill=white] (-1.5,1) circle (.2);
\draw[fill=white] (1.5,-1) circle (.2);
\draw[fill=black] (-1.5,-1) circle (.2);
\draw[fill=white] (1.5,-3) circle (.2);
\draw[fill=white] (-1.5,-3) circle (.2);
\draw[fill=black] (1.5,-5) circle (.2);
\draw[fill=black] (-1.5,-5) circle (.2);

\draw[arrows={-angle 90}, shorten >=2, shorten <=2]  (-1.5,1) -- (1.5,1);
\draw[arrows={-angle 90}, shorten >=2, shorten <=2]  (1.5,1) -- (-1.5,-1);
\draw[arrows={-angle 90}, shorten >=2, shorten <=2]  (-1.5,-1) -- (1.5,-5);

\draw[arrows={-angle 90}, shorten >=2, shorten <=2]  (1.5,-1) -- (-1.5,-3);
\draw[arrows={-angle 90}, shorten >=2, shorten <=2]  (1.5,-3) -- (-1.5,-5);

\end{tikzpicture}
$\quad \quad$
\begin{tikzpicture}[scale=0.3]
\node at (-4,-2) {$p_7 = $};
\draw[fill=black] (1.5,1) circle (.2);
\draw[fill=white] (-1.5,1) circle (.2);
\draw[fill=white] (1.5,-1) circle (.2);
\draw[fill=black] (-1.5,-1) circle (.2);
\draw[fill=white] (1.5,-3) circle (.2);
\draw[fill=black] (-1.5,-3) circle (.2);
\draw[fill=black] (1.5,-5) circle (.2);
\draw[fill=white] (-1.5,-5) circle (.2);

\draw[arrows={-angle 90}, shorten >=2, shorten <=2]  (-1.5,1) -- (1.5,1);
\draw[arrows={-angle 90}, shorten >=2, shorten <=2]  (1.5,1) -- (-1.5,-1);
\draw[arrows={-angle 90}, shorten >=2, shorten <=2]  (-1.5,-1) -- (1.5,-5);

\draw[arrows={-angle 90}, shorten >=2, shorten <=2]  (1.5,-1) -- (-1.5,-3);
\draw[arrows={-angle 90}, shorten >=2, shorten <=2]  (1.5,-3) -- (-1.5,-5);

\end{tikzpicture}
$\quad \quad$
\begin{tikzpicture}[scale=0.3]
\node at (-4,-2) {$p_8 = $};
\draw[fill=black] (1.5,1) circle (.2);
\draw[fill=black] (-1.5,1) circle (.2);
\draw[fill=white] (1.5,-1) circle (.2);
\draw[fill=black] (-1.5,-1) circle (.2);
\draw[fill=white] (1.5,-3) circle (.2);
\draw[fill=white] (-1.5,-3) circle (.2);
\draw[fill=black] (1.5,-5) circle (.2);
\draw[fill=white] (-1.5,-5) circle (.2);

\draw[arrows={-angle 90}, shorten >=2, shorten <=2]  (-1.5,1) -- (1.5,1);
\draw[arrows={-angle 90}, shorten >=2, shorten <=2]  (1.5,1) -- (-1.5,-1);
\draw[arrows={-angle 90}, shorten >=2, shorten <=2]  (-1.5,-1) -- (1.5,-5);

\draw[arrows={-angle 90}, shorten >=2, shorten <=2]  (1.5,-1) -- (-1.5,-3);
\draw[arrows={-angle 90}, shorten >=2, shorten <=2]  (1.5,-3) -- (-1.5,-5);

\end{tikzpicture}
\end{center}
Here the black vertices indicate the corresponding subrepresentation of $M$. The pairs $p_1$ and $p_2$, $p_3$ and $p_4$, and  $p_6$, $p_7$ are each isomorphic as subrepresentations of $M$. The attracting sets of the fixed points are cells by \cite[Theorem~4.13]{Pue2020}. Their dimension equals the number of out going arrows in the following moment graph which is computed using \cite[Theorem~6.15]{LaPu2020}.
\begin{center}
\begin{tikzpicture}[scale=1]
\node at (5,2) {$\mr{with}\ \mr{labels:}$};
\draw[arrows={-angle 90},dash pattern={on 1pt off 2pt on 1pt off 2pt}]  (4.0,1.5) -- (4.5,1.5);
\node at (5.5,1.5) {$\widehat{=} \ \ \epsilon_3-\epsilon_2$};
\draw[arrows={-angle 90},dash pattern={on 4pt off 2pt on4pt off 2pt}]  (4.0,1.0) -- (4.5,1.0);
\node at (5.8,1.0) {$\widehat{=} \ \ \epsilon_2-\epsilon_1-\delta$};
\draw[arrows={-angle 90},dash pattern={on 7pt off 2pt on 1pt off 2pt}]  (4.0,0.5) -- (4.5,0.5);
\node at (5.8,0.5) {$\widehat{=} \ \ \epsilon_3-\epsilon_1-\delta$};
\draw[arrows={-angle 90},dash pattern={on 2pt off 1pt on 1pt off 1pt}]  (4.0,0.0) -- (4.5,0.0);
\node at (5.9,0.0) {$\widehat{=} \ \ \epsilon_1-\epsilon_3+3\delta$};
\draw[arrows={-angle 90},dash pattern={on 2pt off 6pt on 2pt off 6pt}]  (4.0,-0.5) -- (4.5,-0.5);
\node at (5.9,-0.5) {$\widehat{=} \ \ \epsilon_1-\epsilon_2+3\delta$};
\node at (0,2) {$p_7$};
\draw[arrows={-angle 90}, shorten >=8, shorten <=8,dash pattern={on 1pt off 2pt on 1pt off 2pt}]  (0,2) -- (-1.73205,1);
\draw[arrows={-angle 90}, shorten >=8, shorten <=8,dash pattern={on 4pt off 2pt on4pt off 2pt}]  (0,2) -- (1.73205,1);
\draw[arrows={-angle 90}, shorten >=8, shorten <=8,dash pattern={on 7pt off 2pt on 1pt off 2pt}]  (0,2) -- (0,-2);
\node at (-1.73205,3) {$p_8$};
\draw[arrows={-angle 90}, shorten >=8, shorten <=8,dash pattern={on 4pt off 2pt on4pt off 2pt}]  (-1.73205,3) -- (0,2);
\draw[arrows={-angle 90}, shorten >=8, shorten <=8,dash pattern={on 7pt off 2pt on 1pt off 2pt}]  (-1.73205,3) -- (-1.73205,1);
\node at (-1.73205,1) {$p_6$};
\draw[arrows={-angle 90}, shorten >=8, shorten <=8,dash pattern={on 7pt off 2pt on 1pt off 2pt}]  (-1.73205,1) -- (-1.73205,-1);
\draw[arrows={-angle 90}, shorten >=8, shorten <=8,dash pattern={on 4pt off 2pt on4pt off 2pt}]  (-1.73205,1) -- (1.73205,-1);
\node at (1.73205,1) {$p_4$};
\draw[arrows={-angle 90}, shorten >=8, shorten <=8,dash pattern={on 7pt off 2pt on 1pt off 2pt}]  (1.73205,1) -- (1.73205,-1);
\draw[arrows={-angle 90}, shorten >=8, shorten <=8,dash pattern={on 1pt off 2pt on 1pt off 2pt}]  (1.73205,1) -- (-1.73205,-1);
\node at (1.73205,-3) {$p_5$};
\draw[arrows={-angle 90}, shorten >=8, shorten <=8,dash pattern={on 2pt off 1pt on 1pt off 1pt}]  (1.73205,-3) -- (1.73205,-1);
\draw[arrows={-angle 90}, shorten >=8, shorten <=8,dash pattern={on 2pt off 6pt on 2pt off 6pt}]  (1.73205,-3) -- (0,-2);
\node at (-1.73205,-1) {$p_3$};
\draw[arrows={-angle 90}, shorten >=8, shorten <=8,dash pattern={on 4pt off 2pt on4pt off 2pt}]  (-1.73205,-1) -- (0,-2);
\node at (1.73205,-1) {$p_2$};
\draw[arrows={-angle 90}, shorten >=8, shorten <=8,dash pattern={on 1pt off 2pt on 1pt off 2pt}]  (1.73205,-1) -- (0,-2);
\node at (0,-2) {$p_1$};
\end{tikzpicture}
\end{center}
The labels are expressed as linear combination of the characters
\begin{align*}
\epsilon_i &: T \to \C^*; \quad (\gamma_0,\gamma_1,\dots,\gamma_d) \mapsto  \gamma_i \quad \mr{for} \  i \in [d],\\
\delta &: T \to \C^*; \quad (\gamma_0,\gamma_1,\dots,\gamma_d) \mapsto  \gamma_0.
\end{align*}
Here the dashed lines were used to highlight the symmetries of the labeling and avoid to write the labels in the picture.

There are four points which are not rationally smooth. Namely the tangent spaces at $p_1$, $p_2$, $p_6$ and $p_7$ are four-dimensional, whereas $\mr{Gr}_\mb{e}(M)$ itself is three-dimensional. We can read this from the picture as follows: the number of edges adjacent to a point is the dimension of its tangent space and the number of outgoing edges is the cell dimension. The irreducible components are obtained as closure of the strata of the points $p_8$, $p_7$ and $p_5$, because their strata are not contained in the closure of any other stratum. These are generic subrepresentation types of $M$ for dimension vector $\mb{e} = (2,2)$. Hence the desingularization of $\mr{Gr}_\mb{e}(M)$ consists of three components.

The extended quiver $\hat{\Delta}_{2,4}$ is
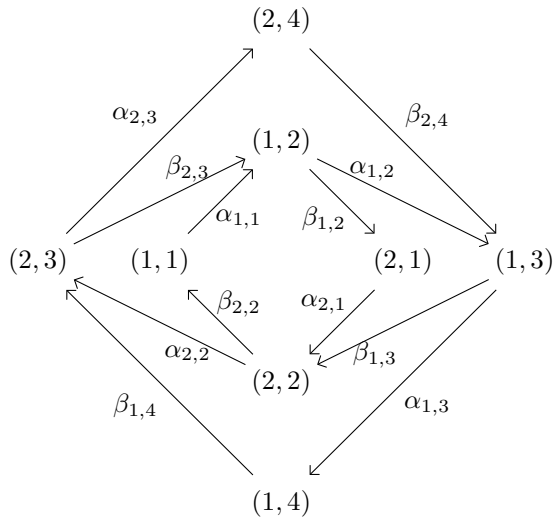
\begin{center}
\begin{tikzpicture}[scale=1.6]
\node at (0,1) {$(1,2)$};
\node at (1,0) {$(2,1)$};
\node at (0,-1) {$(2,2)$};
\node at (-1,0) {$(1,1)$};

\draw[arrows={-angle 90}, shorten >=15, shorten <=15]  (0,1) -- (1,0);\node at (0.35,0.35) {$\beta_{1,2}$};
\draw[arrows={-angle 90}, shorten >=15, shorten <=15]  (1,0) -- (0,-1);\node at (0.35,-0.35) {$\alpha_{2,1}$};
\draw[arrows={-angle 90}, shorten >=15, shorten <=15]  (0,-1) -- (-1,0);\node at (-0.35,-0.35) {$\beta_{2,2}$};
\draw[arrows={-angle 90}, shorten >=15, shorten <=15]  (-1,0) -- (0,1);\node at (-0.35,0.35) {$\alpha_{1,1}$};

\node at (0,2) {$(2,4)$};
\node at (2,0) {$(1,3)$};
\node at (0,-2) {$(1,4)$};
\node at (-2,0) {$(2,3)$};

\draw[arrows={-angle 90}, shorten >=15, shorten <=15]  (0,1) -- (2,0);\node at (0.75,0.75) {$\alpha_{1,2}$};
\draw[arrows={-angle 90}, shorten >=15, shorten <=15]  (0,2) -- (2,0);\node at (1.2,1.2) {$\beta_{2,4}$};

\draw[arrows={-angle 90}, shorten >=15, shorten <=15]  (0,-1) -- (-2,0);\node at (-0.77,-0.77) {$\alpha_{2,2}$};
\draw[arrows={-angle 90}, shorten >=15, shorten <=15]  (0,-2) -- (-2,0);\node at (-1.2,-1.2) {$\beta_{1,4}$};

\draw[arrows={-angle 90}, shorten >=15, shorten <=15]  (-2,0) -- (0,1);\node at (-0.77,0.77) {$\beta_{2,3}$};
\draw[arrows={-angle 90}, shorten >=15, shorten <=15]  (-2,0) -- (0,2);\node at (-1.2,1.2) {$\alpha_{2,3}$};

\draw[arrows={-angle 90}, shorten >=15, shorten <=15]  (2,0) -- (0,-1);\node at (0.77,-0.77) {$\beta_{1,3}$};
\draw[arrows={-angle 90}, shorten >=15, shorten <=15]  (2,0) -- (0,-2);\node at (1.2,-1.2) {$\alpha_{1,3}$};

\end{tikzpicture}
\end{center}
For the basis induced by the basis $B$ of $M$, the coefficient quiver of $\hat{M}$ is
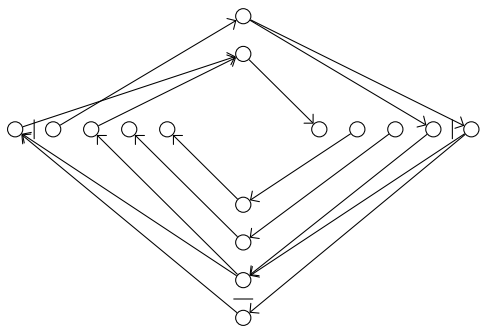
\begin{center}
\begin{tikzpicture}[scale=1]
\draw[fill=white] (3,0) circle (.1);
\draw (2.75,0.13) -- (2.75,-0.13);
\draw[fill=white] (2.5,0) circle (.1);
\draw[fill=white] (2,0) circle (.1);
\draw[fill=white] (1.5,0) circle (.1);
\draw[fill=white] (1,0) circle (.1);

\draw[fill=white] (0,-1) circle (.1);
\draw[fill=white] (0,-1.5) circle (.1);
\draw[fill=white] (0,-2) circle (.1);
\draw (-0.13,-2.25) -- (0.13,-2.25);
\draw[fill=white] (0,-2.5) circle (.1);

\draw[fill=white] (-3,0) circle (.1);
\draw (-2.75,0.13) -- (-2.75,-0.13);
\draw[fill=white] (-2.5,0) circle (.1);
\draw[fill=white] (-2,0) circle (.1);
\draw[fill=white] (-1.5,0) circle (.1);
\draw[fill=white] (-1,0) circle (.1);

\draw[fill=white] (0,1) circle (.1);
\draw[fill=white] (0,1.5) circle (.1);
\draw (-0.13,1.75) -- (0.13,1.75);

\draw[arrows={-angle 90}, shorten >=3, shorten <=3]  (1.5,0) -- (0,-1);
\draw[arrows={-angle 90}, shorten >=3, shorten <=3]  (2,0) -- (0,-1.5);
\draw[arrows={-angle 90}, shorten >=3, shorten <=3]  (2.5,0) -- (0,-2);
\draw[arrows={-angle 90}, shorten >=3, shorten <=3]  (3,0) -- (0,-2);
\draw[arrows={-angle 90}, shorten >=3, shorten <=3]  (3,0) -- (0,-2.5);

\draw[arrows={-angle 90}, shorten >=3, shorten <=3]  (0,-1) -- (-1,0);
\draw[arrows={-angle 90}, shorten >=3, shorten <=3]  (0,-1.5) -- (-1.5,0);
\draw[arrows={-angle 90}, shorten >=3, shorten <=3]  (0,-2) -- (-2,0);
\draw[arrows={-angle 90}, shorten >=3, shorten <=3]  (0,-2) -- (-3,0);
\draw[arrows={-angle 90}, shorten >=3, shorten <=3]  (0,-2.5) -- (-3,0);

\draw[arrows={-angle 90}, shorten >=3, shorten <=3]  (-2,0) -- (0,1);
\draw[arrows={-angle 90}, shorten >=3, shorten <=3]  (-2.5,0) -- (0,1.5);
\draw[arrows={-angle 90}, shorten >=3, shorten <=3]  (-3,0) -- (0,1);

\draw[arrows={-angle 90}, shorten >=3, shorten <=3]  (0,1) -- (1,0);
\draw[arrows={-angle 90}, shorten >=3, shorten <=3]  (0,1.5) -- (2.5,0);
\draw[arrows={-angle 90}, shorten >=3, shorten <=3]  (0,1.5) -- (3,0);

\end{tikzpicture}
\end{center}
Here the separating lines between the vertices indicate if they live over the inner or outer vertex of $\hat{\Delta}_{2,4}$ in that position. Representatives for the extended representations of the generic subrepresentation types are
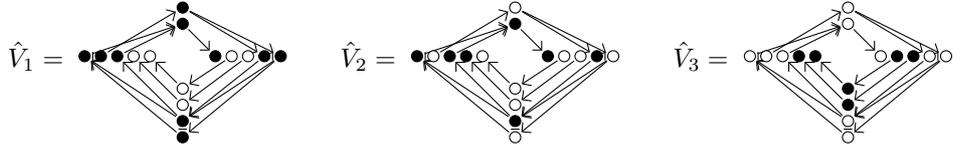
\begin{center}
\begin{tikzpicture}[scale=0.43]
\node at (-4.5,0) {$\hat{V}_1=$};
\draw[arrows={-angle 90}, shorten >=3, shorten <=3]  (1.5,0) -- (0,-1);
\draw[arrows={-angle 90}, shorten >=3, shorten <=3]  (2,0) -- (0,-1.5);
\draw[arrows={-angle 90}, shorten >=3, shorten <=3]  (2.5,0) -- (0,-2);
\draw[arrows={-angle 90}, shorten >=3, shorten <=3]  (3,0) -- (0,-2);
\draw[arrows={-angle 90}, shorten >=3, shorten <=3]  (3,0) -- (0,-2.5);

\draw[arrows={-angle 90}, shorten >=3, shorten <=3]  (0,-1) -- (-1,0);
\draw[arrows={-angle 90}, shorten >=3, shorten <=3]  (0,-1.5) -- (-1.5,0);
\draw[arrows={-angle 90}, shorten >=3, shorten <=3]  (0,-2) -- (-2,0);
\draw[arrows={-angle 90}, shorten >=3, shorten <=3]  (0,-2) -- (-3,0);
\draw[arrows={-angle 90}, shorten >=3, shorten <=3]  (0,-2.5) -- (-3,0);

\draw[arrows={-angle 90}, shorten >=3, shorten <=3]  (-2,0) -- (0,1);
\draw[arrows={-angle 90}, shorten >=3, shorten <=3]  (-2.5,0) -- (0,1.5);
\draw[arrows={-angle 90}, shorten >=3, shorten <=3]  (-3,0) -- (0,1);

\draw[arrows={-angle 90}, shorten >=3, shorten <=3]  (0,1) -- (1,0);
\draw[arrows={-angle 90}, shorten >=3, shorten <=3]  (0,1.5) -- (2.5,0);
\draw[arrows={-angle 90}, shorten >=3, shorten <=3]  (0,1.5) -- (3,0);

\draw[fill=black] (3,0) circle (.17);
\draw (2.75,0.13) -- (2.75,-0.13);
\draw[fill=black] (2.5,0) circle (.17);
\draw[fill=white] (2,0) circle (.17);
\draw[fill=white] (1.5,0) circle (.17);
\draw[fill=black] (1,0) circle (.17);

\draw[fill=white] (0,-1) circle (.17);
\draw[fill=white] (0,-1.5) circle (.17);
\draw[fill=black] (0,-2) circle (.17);
\draw (-0.13,-2.25) -- (0.13,-2.25);
\draw[fill=black] (0,-2.5) circle (.17);

\draw[fill=black] (-3,0) circle (.17);
\draw (-2.75,0.13) -- (-2.75,-0.13);
\draw[fill=black] (-2.5,0) circle (.17);
\draw[fill=black] (-2,0) circle (.17);
\draw[fill=white] (-1.5,0) circle (.17);
\draw[fill=white] (-1,0) circle (.17);

\draw[fill=black] (0,1) circle (.17);
\draw[fill=black] (0,1.5) circle (.17);
\draw (-0.13,1.75) -- (0.13,1.75);
\end{tikzpicture}
$ \quad$
\begin{tikzpicture}[scale=0.43]
\node at (-4.5,0) {$\hat{V}_2=$};
\draw[arrows={-angle 90}, shorten >=3, shorten <=3]  (1.5,0) -- (0,-1);
\draw[arrows={-angle 90}, shorten >=3, shorten <=3]  (2,0) -- (0,-1.5);
\draw[arrows={-angle 90}, shorten >=3, shorten <=3]  (2.5,0) -- (0,-2);
\draw[arrows={-angle 90}, shorten >=3, shorten <=3]  (3,0) -- (0,-2);
\draw[arrows={-angle 90}, shorten >=3, shorten <=3]  (3,0) -- (0,-2.5);

\draw[arrows={-angle 90}, shorten >=3, shorten <=3]  (0,-1) -- (-1,0);
\draw[arrows={-angle 90}, shorten >=3, shorten <=3]  (0,-1.5) -- (-1.5,0);
\draw[arrows={-angle 90}, shorten >=3, shorten <=3]  (0,-2) -- (-2,0);
\draw[arrows={-angle 90}, shorten >=3, shorten <=3]  (0,-2) -- (-3,0);
\draw[arrows={-angle 90}, shorten >=3, shorten <=3]  (0,-2.5) -- (-3,0);

\draw[arrows={-angle 90}, shorten >=3, shorten <=3]  (-2,0) -- (0,1);
\draw[arrows={-angle 90}, shorten >=3, shorten <=3]  (-2.5,0) -- (0,1.5);
\draw[arrows={-angle 90}, shorten >=3, shorten <=3]  (-3,0) -- (0,1);

\draw[arrows={-angle 90}, shorten >=3, shorten <=3]  (0,1) -- (1,0);
\draw[arrows={-angle 90}, shorten >=3, shorten <=3]  (0,1.5) -- (2.5,0);
\draw[arrows={-angle 90}, shorten >=3, shorten <=3]  (0,1.5) -- (3,0);
\draw[fill=white] (3,0) circle (.17);
\draw (2.75,0.13) -- (2.75,-0.13);
\draw[fill=black] (2.5,0) circle (.17);
\draw[fill=white] (2,0) circle (.17);
\draw[fill=white] (1.5,0) circle (.17);
\draw[fill=black] (1,0) circle (.17);

\draw[fill=white] (0,-1) circle (.17);
\draw[fill=white] (0,-1.5) circle (.17);
\draw[fill=black] (0,-2) circle (.17);
\draw (-0.13,-2.25) -- (0.13,-2.25);
\draw[fill=white] (0,-2.5) circle (.17);

\draw[fill=black] (-3,0) circle (.17);
\draw (-2.75,0.13) -- (-2.75,-0.13);
\draw[fill=white] (-2.5,0) circle (.17);
\draw[fill=black] (-2,0) circle (.17);
\draw[fill=black] (-1.5,0) circle (.17);
\draw[fill=white] (-1,0) circle (.17);

\draw[fill=black] (0,1) circle (.17);
\draw[fill=white] (0,1.5) circle (.17);
\draw (-0.13,1.75) -- (0.13,1.75);
\end{tikzpicture}
$ \quad$
\begin{tikzpicture}[scale=0.43]
\node at (-4.5,0) {$\hat{V}_3=$};
\draw[arrows={-angle 90}, shorten >=3, shorten <=3]  (1.5,0) -- (0,-1);
\draw[arrows={-angle 90}, shorten >=3, shorten <=3]  (2,0) -- (0,-1.5);
\draw[arrows={-angle 90}, shorten >=3, shorten <=3]  (2.5,0) -- (0,-2);
\draw[arrows={-angle 90}, shorten >=3, shorten <=3]  (3,0) -- (0,-2);
\draw[arrows={-angle 90}, shorten >=3, shorten <=3]  (3,0) -- (0,-2.5);

\draw[arrows={-angle 90}, shorten >=3, shorten <=3]  (0,-1) -- (-1,0);
\draw[arrows={-angle 90}, shorten >=3, shorten <=3]  (0,-1.5) -- (-1.5,0);
\draw[arrows={-angle 90}, shorten >=3, shorten <=3]  (0,-2) -- (-2,0);
\draw[arrows={-angle 90}, shorten >=3, shorten <=3]  (0,-2) -- (-3,0);
\draw[arrows={-angle 90}, shorten >=3, shorten <=3]  (0,-2.5) -- (-3,0);

\draw[arrows={-angle 90}, shorten >=3, shorten <=3]  (-2,0) -- (0,1);
\draw[arrows={-angle 90}, shorten >=3, shorten <=3]  (-2.5,0) -- (0,1.5);
\draw[arrows={-angle 90}, shorten >=3, shorten <=3]  (-3,0) -- (0,1);

\draw[arrows={-angle 90}, shorten >=3, shorten <=3]  (0,1) -- (1,0);
\draw[arrows={-angle 90}, shorten >=3, shorten <=3]  (0,1.5) -- (2.5,0);
\draw[arrows={-angle 90}, shorten >=3, shorten <=3]  (0,1.5) -- (3,0);
\draw[fill=white] (3,0) circle (.17);
\draw (2.75,0.13) -- (2.75,-0.13);
\draw[fill=white] (2.5,0) circle (.17);
\draw[fill=black] (2,0) circle (.17);
\draw[fill=black] (1.5,0) circle (.17);
\draw[fill=white] (1,0) circle (.17);

\draw[fill=black] (0,-1) circle (.17);
\draw[fill=black] (0,-1.5) circle (.17);
\draw[fill=white] (0,-2) circle (.17);
\draw (-0.13,-2.25) -- (0.13,-2.25);
\draw[fill=white] (0,-2.5) circle (.17);

\draw[fill=white] (-3,0) circle (.17);
\draw (-2.75,0.13) -- (-2.75,-0.13);
\draw[fill=white] (-2.5,0) circle (.17);
\draw[fill=white] (-2,0) circle (.17);
\draw[fill=black] (-1.5,0) circle (.17);
\draw[fill=black] (-1,0) circle (.17);

\draw[fill=white] (0,1) circle (.17);
\draw[fill=white] (0,1.5) circle (.17);
\draw (-0.13,1.75) -- (0.13,1.75);
\end{tikzpicture}
\end{center}
With the explicit description of the cellular decomposition of the quiver Grassmannians $\mr{Gr}_{\bdim \hat{V}_1}(\hat{M})$, $\mr{Gr}_{\bdim \hat{V}_2}(\hat{M})$ and $\mr{Gr}_{\bdim \hat{V}_3}(\hat{M})$ from Theorem~\ref{trm:cellular-decomp-desing}, it is a straightforward computation that their moment graphs are
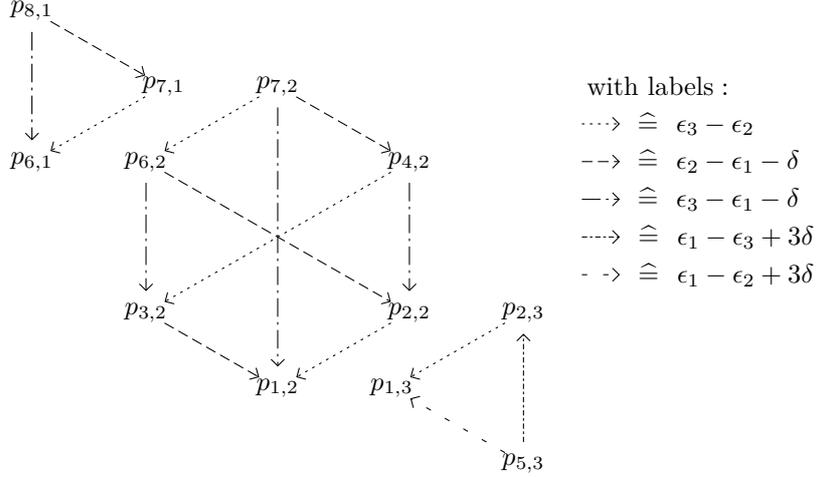
\begin{center}
\begin{tikzpicture}[scale=1]
\node at (5,2) {$\mr{with}\ \mr{labels:}$};
\draw[arrows={-angle 90},dash pattern={on 1pt off 2pt on 1pt off 2pt}]  (4.0,1.5) -- (4.5,1.5);
\node at (5.5,1.5) {$\widehat{=} \ \ \epsilon_3-\epsilon_2$};
\draw[arrows={-angle 90},dash pattern={on 4pt off 2pt on4pt off 2pt}]  (4.0,1.0) -- (4.5,1.0);
\node at (5.8,1.0) {$\widehat{=} \ \ \epsilon_2-\epsilon_1-\delta$};
\draw[arrows={-angle 90},dash pattern={on 7pt off 2pt on 1pt off 2pt}]  (4.0,0.5) -- (4.5,0.5);
\node at (5.8,0.5) {$\widehat{=} \ \ \epsilon_3-\epsilon_1-\delta$};
\draw[arrows={-angle 90},dash pattern={on 2pt off 1pt on 1pt off 1pt}]  (4.0,0.0) -- (4.5,0.0);
\node at (5.9,0.0) {$\widehat{=} \ \ \epsilon_1-\epsilon_3+3\delta$};
\draw[arrows={-angle 90},dash pattern={on 2pt off 6pt on 2pt off 6pt}]  (4.0,-0.5) -- (4.5,-0.5);
\node at (5.9,-0.5) {$\widehat{=} \ \ \epsilon_1-\epsilon_2+3\delta$};
\node at (0-1.5,2) {$p_{7,1}$};
\draw[arrows={-angle 90}, shorten >=8, shorten <=8,dash pattern={on 1pt off 2pt on 1pt off 2pt}]  (0-1.5,2) -- (-1.73205-1.5,1);
\node at (-1.73205-1.5,3) {$p_{8,1}$};
\draw[arrows={-angle 90}, shorten >=8, shorten <=8,dash pattern={on 4pt off 2pt on4pt off 2pt}]  (-1.73205-1.5,3) -- (0-1.5,2);
\draw[arrows={-angle 90}, shorten >=8, shorten <=8,dash pattern={on 7pt off 2pt on 1pt off 2pt}]  (-1.73205-1.5,3) -- (-1.73205-1.5,1);
\node at (-1.73205-1.5,1) {$p_{6,1}$};

\node at (0,2) {$p_{7,2}$};
\draw[arrows={-angle 90}, shorten >=8, shorten <=8,dash pattern={on 1pt off 2pt on 1pt off 2pt}]  (0,2) -- (-1.73205,1);
\draw[arrows={-angle 90}, shorten >=8, shorten <=8,dash pattern={on 4pt off 2pt on4pt off 2pt}]  (0,2) -- (1.73205,1);
\draw[arrows={-angle 90}, shorten >=8, shorten <=8,dash pattern={on 7pt off 2pt on 1pt off 2pt}]  (0,2) -- (0,-2);
\node at (-1.73205,1) {$p_{6,2}$};
\draw[arrows={-angle 90}, shorten >=8, shorten <=8,dash pattern={on 7pt off 2pt on 1pt off 2pt}]  (-1.73205,1) -- (-1.73205,-1);
\draw[arrows={-angle 90}, shorten >=8, shorten <=8,dash pattern={on 4pt off 2pt on4pt off 2pt}]  (-1.73205,1) -- (1.73205,-1);
\node at (1.73205,1) {$p_{4,2}$};
\draw[arrows={-angle 90}, shorten >=8, shorten <=8,dash pattern={on 7pt off 2pt on 1pt off 2pt}]  (1.73205,1) -- (1.73205,-1);
\draw[arrows={-angle 90}, shorten >=8, shorten <=8,dash pattern={on 1pt off 2pt on 1pt off 2pt}]  (1.73205,1) -- (-1.73205,-1);
\node at (-1.73205,-1) {$p_{3,2}$};
\draw[arrows={-angle 90}, shorten >=8, shorten <=8,dash pattern={on 4pt off 2pt on4pt off 2pt}]  (-1.73205,-1) -- (0,-2);
\node at (1.73205,-1) {$p_{2,2}$};
\draw[arrows={-angle 90}, shorten >=8, shorten <=8,dash pattern={on 1pt off 2pt on 1pt off 2pt}]  (1.73205,-1) -- (0,-2);
\node at (0,-2) {$p_{1,2}$};

\node at (2.73205+0.5,-3) {$p_{5,3}$};
\draw[arrows={-angle 90}, shorten >=8, shorten <=8,dash pattern={on 2pt off 1pt on 1pt off 1pt}]  (2.73205+0.5,-3) -- (2.73205+0.5,-1);
\draw[arrows={-angle 90}, shorten >=8, shorten <=8,dash pattern={on 2pt off 6pt on 2pt off 6pt}]  (2.73205+0.5,-3) -- (1+0.5,-2);
\node at (2.73205+0.5,-1) {$p_{2,3}$};
\draw[arrows={-angle 90}, shorten >=8, shorten <=8,dash pattern={on 1pt off 2pt on 1pt off 2pt}]  (2.73205+0.5,-1) -- (1+0.5,-2);
\node at (1+0.5,-2) {$p_{1,3}$};
\end{tikzpicture}
\end{center}
Here $\hat{p}_{i,j}$ is the preimage of $p_i$ in $\mr{Gr}_{\bdim \hat{V}_j}(\hat{M})$. Moreover from the cellular decompositions we obtain the isomorphisms 
\begin{align*}
    \mr{Gr}_{\bdim \hat{V}_1}(\hat{M}) &\cong \mr{Gr}_1(\C^3),\\
    \mr{Gr}_{\bdim \hat{V}_2}(\hat{M}) &\cong \mathcal{F}l(\mr{SL}_3),\\
    \mr{Gr}_{\bdim \hat{V}_3}(\hat{M}) &\cong \mr{Gr}_2(\C^3).
\end{align*}
With the moment graph of the desingularization as described above, it is possible to compute the Euler classes at the singular points of $\mr{Gr}_\mb{e}(M)$ using Lemma~\ref{lma:euler-class-along-resolution}. For example we obtain
\begin{align*}
    \mr{Eu}_T(p_1,Z_5) &= \frac{1}{(\epsilon_3-\epsilon_2)(\epsilon_2-\epsilon_1-\delta)}+ \frac{1}{(\epsilon_3-\epsilon_2)(\epsilon_1-\epsilon_2+3\delta)}\\
    &=\frac{2\delta}{(\epsilon_3-\epsilon_2)(\epsilon_2-\epsilon_1-\delta)(\epsilon_1-\epsilon_2+3\delta)}
\end{align*}
where 
\( Z_5 = \cup_{i=1}^5W_i\).

We compute the following basis of $H_T^\bullet(\mr{Gr}_\mb{e}(M))$ as free module over $H_T^\bullet(pt)$:
\begin{align*}
    \varphi^{(1)}&=(1,1,1,1,1,1,1,1)\\
    \varphi^{(2)}&=(0,\epsilon_3-\epsilon_2,0,\epsilon_3-\epsilon_2,\epsilon_1-\epsilon_2+3\delta,\epsilon_3-\epsilon_1-\delta,\epsilon_3-\epsilon_1-\delta,\epsilon_3-\epsilon_1-\delta)\\
    \varphi^{(3)}&=(0,0,\epsilon_2-\epsilon_1-\delta,\epsilon_3-\epsilon_1-\delta,0,\epsilon_2-\epsilon_1-\delta,\epsilon_3-\epsilon_2,\epsilon_3-\epsilon_1-\delta)\\
    \varphi^{(4)}&=(\epsilon_3-\epsilon_2)(\epsilon_3-\epsilon_1-\delta)\cdot(0,0,0,1,0,0,1,0)\\
    \varphi^{(5)}&=(\epsilon_1-\epsilon_2+3\delta)(\epsilon_1-\epsilon_3+3\delta)\cdot(0,0,0,0,1,0,0,0)\\
    \varphi^{(6)}&=(\epsilon_2-\epsilon_1-\delta)(\epsilon_3-\epsilon_1-\delta)\cdot(0,0,0,0,0,1,1,0)\\
    \varphi^{(7)}&=(\epsilon_3-\epsilon_2)(\epsilon_2-\epsilon_1-\delta)(\epsilon_3-\epsilon_1-\delta)\cdot(0,0,0,0,0,0,1,0)\\
    \varphi^{(8)}&=(\epsilon_2-\epsilon_1-\delta)(\epsilon_3-\epsilon_1-\delta)\cdot(0,0,0,0,0,0,0,1)
\end{align*}
Observe that the special role of $p_8$ in this example allows to generate more zero-entries as in the general setting of Theorem~\ref{trm:cohomology-generators-general-setting}.

\end{document}